\documentclass[12pt]{article}
\usepackage{textcomp}
\usepackage{graphicx}           
\usepackage{epstopdf}
\usepackage{color}              
\usepackage{xcolor,graphicx}
\usepackage{geometry}
\usepackage{float}
\usepackage{indentfirst}        
\usepackage{amsmath,amssymb,bm} 
\usepackage{amsthm}
\usepackage{cases}              
\usepackage{setspace}
\usepackage{hyperref}
\usepackage{titlesec}
\usepackage[all]{xy}            
\usepackage{tikz}
\usepackage{multirow}
\usepgflibrary{arrows}
\hypersetup{dvips,ps2pdf,CJKbookmarks,hypertex,
	colorlinks=true,
	pdfpagemode=FullScreen,
	bookmarksopen=true,
	bookmarksnumbered=true}     
\usepackage{graphicx} 
\usepackage{epstopdf}
\usepackage{caption}
\usepackage{diagbox}
\usepackage{mathrsfs}

\usepackage{dsfont}
\usepackage{graphics}
\usepackage{hyperref,cleveref,cite}
\usepackage{titlesec}
\usepackage{tcolorbox}
\usepackage{fancyhdr,fancyvrb}
\usepackage{xcolor,graphicx,geometry,float,indentfirst,cases,setspace,titlesec}
\usepackage{diagbox}
\usepackage{times}

\newtheorem{con}{Conjecture}[section]

\newtheorem{lem}[con]{Lemma}
\newtheorem{prop}[con]{Proposition}
\newtheorem{coro}[con]{Corollary}
\newtheorem{thm}[con]{Theorem}

\newtheorem{pf}{proof}[section]
\theoremstyle{remark}

\begingroup
\makeatletter \@addtoreset{equation}{section} \makeatother
\makeindex 
\def\qed{\hfill \rule{4pt}{7pt}}
\def\pf{\vskip 0.2cm {\noindent \bf Proof.}\quad}

\begin{document}
\begin{center}

 {\Large \bf Signed Mahonian Polynomials on Derangements \\[5pt] in  Classical Weyl Groups}

\end{center}

\begin{center}
{Kathy Q. Ji}$^{1}$ and {Dax T.X. Zhang}$^{2}$   \vskip 2mm

$^1$ Center for Applied Mathematics\\[3pt]
Tianjin University\\[3pt]
Tianjin 300072, P.R. China\\[6pt]
  kathyji@tju.edu.cn
   \vskip 2mm
$^2$ College of Mathematical Science \\[3pt]
Institute of Mathematics and Interdisciplinary Sciences\\[3pt]
Tianjin Normal University\\[3pt]
Tianjin 300387, P.R. China\\[6pt]
   \vskip 2mm
   zhangtianxing6@tju.edu.cn
\end{center}
\vskip 2mm
\begin{center}
{\bf Abstract}
\end{center}
The polynomial of the   major index ${\rm maj}_W (\sigma)$  over the subset $T$ of the Coxeter group $W$ is called the Mahonian polynomial over $T$,   where ${\rm maj}_W (\sigma)$ is a   Mahonian statistic of an element $\sigma \in T$,  whereas the polynomial of the   major index ${\rm maj}_W (\sigma)$ with the sign $(-1)^{\ell_W(\sigma)}$  over the subset $T$    is referred to as  the signed  Mahonian polynomial over $T$, where ${\ell_W(\sigma)}$ is the length   of  $\sigma \in T$.  Gessel, Wachs, and Chow established formulas for the Mahonian polynomials over the sets of derangements in the symmetric group $S_n$ and the hyperoctahedral group $B_n$. By extending Wachs' approach and employing a refinement of Stanley's shuffle theorem established in our recent paper\cite{Ji-Zhang-2024}, we derive a formula for the Mahonian polynomials over the set of derangements in the even-signed permutation group $D_n$.  This completes a picture which is now known for all the classical Weyl groups.  Gessel-Simion, Adin-Gessel-Roichman, and Biagioli previously established formulas for the signed Mahonian polynomials over the classical Weyl groups.  Building upon their formulas,
we derive some new formulas for the signed Mahonian polynomials  over the set of derangements in classical Weyl groups. As applications of the formulas for the (signed) Mahonian polynomials    over  the sets of derangements in the classical Weyl groups, we obtain enumerative formulas of the number of derangements in classical Weyl groups with even lengths.

\vskip 2mm

\noindent
{\bf Keywords:}  permutations,  signed permutations,  derangement, length functions, Mahonian statistic, signed Mahonian polynomials, Classical Weyl groups

\noindent
{\bf AMS Classification:}   05A05, 05A15, 05A19, 05A30

 \vskip 2mm

\section{Introduction}

This paper is concerned with the $q$-counting derangements in classical Weyl groups  by their major indices.   The classical Weyl groups may be described as follows: $S_n$ is the symmetric group consisting of all permutations of the set $[n]:=\{1,2,\ldots, n\}$, and $B_n$ is the hyperoctaherdral
 group consisting of all signed permutations of   $[n]$, and $D_n$ is the subgroup of index two in $B_n$ consisting of signed permutations of   $[n]$ with an even number of signs.
 Recall that a permutation of $[n]$ is a bijection $\sigma \colon [n]\rightarrow [n]$,  whereas a signed permutation of $[n]$ is defined to be a function $\sigma \colon [n]\rightarrow [-n,n]\setminus \{0\}$ such that $|\sigma|$ is a permutation of $[n]$, where $|\sigma|(i)=|\sigma(i)|$ for $i \in [n]$.  More precisely,   an element $\sigma$ of $S_n$ can be presented as $\sigma=\sigma_1\sigma_2\cdots \sigma_n$, where $\sigma_i=\sigma(i)$ and the an element $\sigma$ of $B_n$ (or $D_n$) as $\sigma=\sigma_1\sigma_2\cdots \sigma_n$, where some elements are associated with the minus sign.  For convenience, we write $\overline{i}=-i$.

 We can also regard each of the classical Weyl groups $S_n$, $B_n$ and $D_n$ {as a Coxeter group $W$ with a generating set $S$}, see Bj\"orner and Brenti \cite{Bjorner-Brenti-2005} for background. Each element $w\in W$ can be written as a product of generators: $w=s_1s_2\cdots s_k$ where $s_i \in S.$ 
 
 For an element  $\sigma \in W$, we define its length $\ell_W(\sigma)$ by
\begin{equation}
    \ell_W(\sigma)= \min\{k \colon \sigma = s_{i_1}\cdots s_{i_k} \quad \text{for some } s_{i_j}\in S \}.
\end{equation}
The descent set of $\sigma$, denoted ${\rm Des}_W(\sigma)$,  is defined by
\begin{equation}
    {\rm Des}_W(\sigma)= \{i \colon \ell_W(\sigma s_i)<\ell_W(\sigma)\}.
\end{equation}
And the number of  descents of $\sigma$, denoted ${\rm des}_W(\sigma)$,  is given by
\begin{equation}
    {\rm des}_W(\sigma)= \#{\rm Des}_W(\sigma).
\end{equation}
Let ${\rm maj}_W(\sigma)$ denote a Mahonian statistic of $\sigma$, which is equidistributed with the length  $\ell_W(\sigma)$. To wit,
\begin{equation}\label{defi-poincare}
\sum_{\sigma \in W}    q^{{\rm maj}_W(\sigma)}=\sum_{\sigma \in W}    q^{{\rm \ell}_W(\sigma)}.
\end{equation}
The  polynomial \eqref{defi-poincare} is known as the Poincare\'e polynomial of the group $W$, which has a nice product formula for every finite Coxeter group, see Bj\"orner and Brenti \cite[Chapter 7]{Bjorner-Brenti-2005}.

The signed Mahonian polynomial over the Coxeter group $W$ is defined as
\begin{equation} \label{main-pol}
\sum_{\sigma \in W}   (-1)^{\ell
_W(\sigma)} q^{{\rm maj}_W(\sigma)}.
\end{equation}
 Gessel and Simion (see \cite[Corollay 2]{Wachs-1992}) first obtained an elegant factorial-type formula for the signed Mahonian polynomial over  the symmetric group $S_n$ (see Theorem \ref{S-1-thm}  below).    Adin, Gessel and Roichman \cite{Adin-Gessel-Roichman-2005} derived the formula for the signed Mahonian polynomial over the hyperoctaherdral group $B_n$  (see Theorem \ref{BMS-1} below),
and the formula for the signed Mahonian polynomial over the even-signed permutation group  $D_n$ was established by Biagioli \cite{Biagioli-2006}  (see Theorem \ref{sigMahothm} below).

It is worth noting that the enumeration of the descent number ${\rm des}_W(\sigma)$ of the Coxeter group $W$ with the sign $(-1)^{\ell_W(\sigma)}$,  known as the signed Eulerian polynomials,  has been explored by
D\'esarm\'enien-Foata \cite{Desarmenien-Foata-1992}, Loday \cite{Loday-1989},  Reiner \cite{Reiner-1995}  and Wachs \cite{Wachs-1992}. Recently,  the study of the signed Mahonian polynomials over other groups has been  undertaken by Biagioli-Caselli \cite{Biagioli-Caselli-2012},  Caselli \cite{Caselli-2012},  Chang-Eu-Fu-Lin-Lo \cite{Chang-Eu-Fu-Lin-Lo-2019}, Eu-Fu-Hsu-Liao-Sun \cite{Eu-Fu-Hsu-Liao-Sun-2020}  and Eu-Fu-Hsu-Lo \cite{Eu-Fu-Hsu-Lo-2022}.

The main objective of this work is to investigate the polynomial \eqref{main-pol} defined on the set of  the derangements in  classical Weyl groups. Let $\sigma=\sigma_1\cdots \sigma_n \in W$. We say that $i$ is  a fixed point of $\sigma$ if $\sigma_i=i$.  A derangement of type $A$ (resp. $B$ or $D$) is a permutation $\sigma=\sigma_1\cdots \sigma_n$ in $S_n$ (resp. $B_n$ or $D_n$) such that $\sigma_i\neq i$ for $i \in [n]$. The derangement of type $B$ or $D$ is also called the signed derangement.  Denote by $\mathcal{D}_n^S$,  $\mathcal{D}^{B}_n$ and $\mathcal{D}^{D}_n$ the sets of derangements in $S_n$, $B_n$ and $D_n$ respectively.   For example,
\[\mathcal{D}^S_2=\{2\, 1\}, \quad \mathcal{D}^B_2=\{\bar{1}\,\bar{2},2\,1,2\,\bar{1},\bar{2}\,1,\bar{2}\,\bar{1}\}, \quad \text{and} \quad \mathcal{D}^{D}_2=\{\bar{1}\,\bar{2},2\,1,\bar{2}\,\bar{1}\}.\]
The enumeration of derangements has a rich history with notable developments  in the study of derangements by specific permutation statistics. Significant contributions in this direction include the work of Brenti \cite{Brenti-1994}, Gessel \cite{Gessel-Reutenauer-1993} and Wachs \cite{Wachs-1989}, who enumerated derangements by their weak excedances and major indices. It is noteworthy that the $q$-counting of derangements based on their weak excedances has been extensively explored by Assaf \cite{Assaf-2010}, Brenti \cite{Brenti-1994}, Chow \cite{Chow-2009}, Chen-Tang-Zhao \cite{Chen-Tang-Zhao-2009}, Chow-Mansour \cite{Chow-Mansour-2010}, and Pei-Zeng \cite{Pei-Zeng-2024}.

In this paper, our primary focus is on the $q$-counting of derangements in classical Weyl groups based on their major indices. Specifically, we investigate the Mahonian polynomials defined on the set $\mathcal{D}^W_n$ of derangements in the classical Weyl group $W$:
 \begin{equation} \label{main-pol-derang}
\sum_{\sigma \in \mathcal{D}_n^W}   q^{{\rm maj}_W(\sigma)}.
\end{equation}
Gessel (published in \cite{Gessel-Reutenauer-1993}) first obtained an elegant formula for the polynomial \eqref{main-pol-derang} when $W=S_n$ (see Theorem \ref {dSn-Wachs} below), which has been studied further by Chen-Rota \cite{Chen-Rota-1992}, Garsia-Remmel \cite{Garsia-Remmel-1980} and  Wachs \cite{Wachs-1989}. In particular,  Wachs \cite{Wachs-1989} found a bijection on $S_n$ by rearranging a
permutation $\pi \in S_n$ according to excedant ($\pi_i>i$), fixed point, and subcedant ($\pi_i<i$). Wachs demonstrated that this bijection preserves the major index. Subsequently, Garsia and Gessel's result on shuffles of permutations was applied to derive the formula for the polynomial \eqref{main-pol-derang} over $\mathcal{D}^S_n$.

Chow \cite{Chow-2006}  generalized Wachs' approach to derive a formula for the polynomial \eqref{main-pol-derang} when $W=B_n$  (see Theorem \ref{dBn} below). Recently, Chow \cite{Chow-2023}  explored the enumeration of derangements in $D_n$  based on a result of Foata-Han  \cite{Foata-Han-2008}, but Chow's $q$-counting of   derangements in $D_n$   relies on  the natural maj index on $D_n$, not on the Mahonian major index ${\rm maj}_D$  (see Theorem \ref {ThDnchow} below). In this paper, we obtain a formula for the polynomial \eqref{main-pol-derang} when $W=D_n$  by generalizing Wachs' approach. It should be stressed that a refinement of Garsia-Gessel's formula on shuffles of permutations is necessary in the establishment of the formula for the polynomial \eqref{main-pol-derang} when $W=D_n$. This refinement is derived from a refinement of Stanley's shuffle theorem established in our recent paper \cite{Ji-Zhang-2024}.

By utilizing Gessel-Simion,  Adin-Gessel-Roichman and Biagioli's formulas for the signed Mahonian polynomial \eqref{main-pol} over the Coxeter group $W$,  we derive the formulas for the following signed Mahonian polynomials defined on   $\mathcal{D}^W_n$,
\begin{equation} \label{main-pol-derangs}
\sum_{\sigma \in \mathcal{D}_n^W}   (-1)^{\ell
_W(\sigma)} q^{{\rm maj}_W(\sigma)}.
\end{equation}
See Theorem \ref {ThSn}, Theorem \ref{ThBn} and Theorem  \ref{ThDn-} below.

Combining the formulas for the   Mahonian polynomial \eqref{main-pol-derang}   and the signed Mahonian polynomial \eqref{main-pol-derangs}, one leads to the formulas for  the   Mahonian polynomials defined on  $\mathcal{D}^{AS}_n$,  $\mathcal{D}^{AB}_n$ and $\mathcal{D}^{AD}_n$, where   $\mathcal{D}^{AS}_n$,  $\mathcal{D}^{AB}_n$ and $\mathcal{D}^{AD}_n$ denotes the sets of derangements with even length in $S_n$, $B_n$ and $D_n$ respectively (see Corollary \ref{ThAAn}, Corollary \ref{ThABn} and Corollary \ref{ThADn} below).     For example,
\[\mathcal{D}^{AS}_2=\emptyset, \quad \mathcal{D}^{AB}_2=\{\bar{1}\,\bar{2},2\,\bar{1},\bar{2}\,1\}, \quad {\rm and}\quad \mathcal{D}^{AD}_2=\{\bar{1}\,\bar{2}\}.\]
As applications of these   formulas for the Mahonian polynomials over  $\mathcal{D}^{AS}_n$,  $\mathcal{D}^{AB}_n$ and $\mathcal{D}^{AD}_n$, we derive the following enumerative results involving $\mathcal{D}^{AS}_n$,  $\mathcal{D}^{AB}_n$ and $\mathcal{D}^{AD}_n$. It would be interesting to provide combinatorial proofs of these enumerative results.

\begin{thm}\label{en-as} Let $d^{AS}_n$ denote the number of derangements with even length in $S_n$.

{\rm (i)} For $n\geq 2$,

\begin{equation}
d^{AS}_n =\frac{n!}{2}\sum_{k=0}^{n-2}\frac{(-1)^k}{k!}+(-1)^{n-1}(n-1);
\end{equation}

{\rm (ii)} For $n\geq 2$,
\begin{equation}
d^{AS}_n=n d^{AS}_{n-1}+\frac{(-1)^{n-1}}{2}(n-2)(n+1)
\end{equation}
with initial condition $d^{AS}_1=0${\rm ;}

{\rm (iii)} For $n\geq 3$,
\begin{equation}
d^{AS}_n=(n-1)(d^{AS}_{n-1}+d^{AS}_{n-2}+(-1)^{n-1})
\end{equation}
with initial conditions $d^{AS}_1=0$ and $d^{AS}_2=0$.
\end{thm}

 \begin{thm} \label{en-ab} Let $d^{AB}_n$ denote the number of derangements with even length in $B_n$.

{\rm (i)} For $n\geq 1$,
\begin{equation}
d^{AB}_n=n!\sum_{k=0}^{n-1}\frac{2^{n-k-1}(-1)^k}{k!}+(-1)^{n}{\rm ;}
\end{equation}

{\rm (ii)} For $n\geq 2$,
\begin{equation}
d^{AB}_n=2n d^{AB}_{n-1}+(-1)^n(n+1)
\end{equation}
with initial condition $d^{AB}_1=0${\rm ;}

{\rm (iii)} For $n\geq 3$,
\begin{equation}
d^{AB}_n=(n-1)(2d^{AB}_{n-1}+4d^{AB}_{n-2}+(-1)^{n-1})
\end{equation}
 with initial conditions $d^{AB}_1=0$ and $d^{AB}_2=3$.
\end{thm}

\begin{thm} \label{en-ad} Let $d^{AD}_n$ denote the number of derangements with even length in $D_n$.

{\rm (i)} For $n\geq 2$,
\begin{equation}
d^{AD}_n=n!\sum_{k=0}^{n-2}\frac{2^{n-k-2}(-1)^k}{k!}+(-1)^{n-1}(n-1){\rm ;}
\end{equation}

{\rm (ii)} For $n\geq 2$,
\begin{equation}
d^{AD}_n=2n d^{AD}_{n-1}+(-1)^{n-1}(n^2-2n-1)
\end{equation}
with initial condition $d^{AD}_1=0${\rm ;}

{\rm (iii)} For $n\geq 3$,
\begin{align}
d^{AD}_n
&=(2n-1) d^{AD}_{n-1}+2(n-1)d^{AD}_{n-2}+(-1)^{n-1}(2n-3).
\end{align}
with initial conditions $d^{AD}_1=0$ and $d^{AD}_2=1$
\end{thm}

The rest of the paper is organized as follows. Section 2 is dedicated to exploring  the signed Mahonian polynomials over the set of derangements in $S_n$ with the aid of the formula for the signed Mahonian polynomials over $S_n$ established by Gessel-Simion. In Section 3, we begin  by revisiting the combinatorial definitions of the length  $\ell_B(\sigma)$ and the Mahonian major index ${\rm maj}_B(\sigma)$ of an element $\sigma$ in the hyperoctahedral group $B_n$. Subsequently, we derive the formula for signed Mahonian polynomials over the set of  derangements of type $B$ building upon the formula for signed Mahonian polynomials over $B_n$ established by Adin-Gessel-Roichman.  In Section 4, we recall the combinatorial definitions of the length   $\ell_D(\sigma)$ and the Mahonian major index ${\rm maj}_D(\sigma)$ introduced by Biagioli-Caselli for the element $\sigma$ in the even-signed permutation group $D_n.$ We then state the formulas for the Mahonian polynomials and the signed Mahonian polynomials over the set of  signed  derangements   in $D_n$. The establishment of these formulas for these two polynomials over the set of signed  derangements   in $D_n$ is shown to reduce to the derivation of the formulas for   the (signed) Mahonian polynomials over the set of signed derangements in $\Delta_n=\{\gamma=\gamma_1\cdots \gamma_n \in B _ { n }, \gamma_n>0\}$, as detailed in Section 5.  It should be stressed that a refinement of Stanley's shuffle theorem, obtained in our recent paper  \cite{Ji-Zhang-2024}, plays a crucial role in deriving the formula for the (signed)  Mahonian polynomials  over the set of signed derangements in $\Delta_n.$

\section{The signed Mahonian polynomials over $\mathcal{D}_n^S$}

This section is devoted to investigating the signed Mahonian polynomials over the set of derangements in the symmetric group $S_n$. The symmetric group $S_n$ can be viewed as a Coxeter group of type $A$ with the set of generators $S = \{s_1 ,\ldots,s_{n-1} \}$, where $s_i=[1,2,\ldots,i-1,i+1,i,i+2,\ldots,n]$, see \cite[Proposition 1.5.4]{Bjorner-Brenti-2005}. 

Let $\sigma=\sigma_1\cdots \sigma_n \in S_n$. Define
\begin{eqnarray}
 {\rm inv}(\sigma)&=& \{(i,j) \colon i<j \quad \text{and} \quad \sigma_i>\sigma_{j}\},\label{def-inv} \label{defi-inv}\\[5pt]
    {\rm Des}(\sigma)&=& \{1\leq i\leq n-1 \colon \sigma_i>\sigma_{i+1}\}, \label{def-Des} \\[5pt]
    {\rm des}(\sigma)&=& \# {\rm Des}(\sigma), \label{def-des}\\[5pt]
    {\rm maj}(\sigma)&=&\sum_{i \in {\rm Des}(\sigma)} i. \label{def-maj}
\end{eqnarray}
For example, for the permutation $\sigma=5\,3\,1\,2\,4 \in S_5$, we have
\[{\rm inv}(\sigma)=6, \quad {\rm des}(\sigma)=2, \quad \text{and} \quad  {\rm maj}(\sigma)=3.\]
  The length function $\ell_A(\sigma)$  and the descent set ${\rm Des}_A(\sigma)$ of $\sigma \in S_n$ can be computed directly as follows:
 \begin{equation}
 \ell_A(\sigma)={\rm inv}(\sigma),  \quad \text{and} \quad  {\rm Des}_A(\sigma)={\rm Des}(\sigma).
 \end{equation}
A well-known result due to MacMahon \cite{MacMahon-1915} asserts that the major index defined as \eqref{def-maj}  is a Mahonian statistic of $\sigma$. Namely,
\begin{thm}[MacMahon]\label{MacMahon} For $n\geq 1$,
\begin{equation}
    \sum_{\sigma \in S_n}q^{{\rm maj}(\sigma)}=\sum_{\sigma \in S_n}q^{{\rm inv}(\sigma)}=\sum_{\sigma \in S_n}q^{\ell_A(\sigma)}=[n]_q!.
\end{equation}
\end{thm}
Here and in the sequel,  for  a positive integer $n$, we define
\[[n]_q:=\frac{1-q^n}{1-q}=1+q+\cdots +q^{n-1}\]
and for $n\geq 1$, 
\[[n]_q!:=[1]_q[2]_q\cdots [n]_q.\]
{Assume that $[0]_q!=1$.}

Gessel and Simion \cite[Corollary 2]{Wachs-1992} were the first to investigate the signed Mahonian polynomial over the symmetric group $S_n$ and derived the following elegant factorial-type product formula:
\begin{thm}[Gessel-Simion]  \label{S-1-thm} For $n\geq 1$,
\begin{align}\label{S-1}
    \sum_{ \sigma \in S_n} (-1)^{\ell_A(\sigma)} q^{{\rm maj}(\sigma)}&= [1]_q [2]_{-q}[3]_q [4]_{-q} \cdots [n]_{(-1)^{n-1}q} \nonumber \\[5pt]
    &=\left(\frac{1-q} {1+q}\right)^{\lfloor \frac{n}{2}\rfloor}[n]_q!.
\end{align}
\end{thm}

The elegant formula for the Mahonian polynomial over $\mathcal{D}_n^S$ was initially derived by Gessel and Reutenauer published in \cite{Gessel-Reutenauer-1993} as a consequence of the quasi-symmetric generating function encoding the descents and the cycle structure of permutations. Wachs \cite{Wachs-1989} later provided a beautiful combinatorial proof of this formula, utilizing a shuffle theorem due to Garsia and Gessel \cite{Garsia-Gessel-1979}.
\begin{thm}[Gessel-Reutenauer-Wachs]\label{dSn-Wachs} For $n\geq 1$,
\begin{equation}\label{dSn-Wachs-eq}
    d^{S}_n(q)=\sum_{\sigma \in \mathcal{D}_n^S}q^{{\rm maj}(\sigma)}=[n]_q!\sum^n_{k=0}\frac{(-1)^k}{[k]_q!}q^{\binom{k}{2}}.
\end{equation}
\end{thm}

The first result of this paper is to establish the following formula for the signed  Mahonian polynomial over $\mathcal{D}_n^S$.
\begin{thm}\label{ThSn}For $n\geq 1$,
    \begin{equation*}
   \overline{d}^{\,S}_n(q)=\sum_{\sigma \in \mathcal{D}_n^S} (-1)^{\ell_A(\sigma)} q^{{\rm maj}(\sigma)}=[n]_q!\sum^n_{k=0}\frac{(-1)^k}{[k]_q!}q^{\binom{k}{2}}\left(\frac{1-q} {1+q}\right)^{\lfloor \frac{n-k}{2}\rfloor}.
\end{equation*}
\end{thm}
To  establish Theorem \ref{ThSn},  let us first review the combinatorial settings of Wachs.  Let $A=\{0<a_1<a_2<\cdots <a_n\}$ and let $\mathfrak{S}_A$ denote the set of permutations of the set $A$. For
 $\pi=\pi_1\cdots \pi_n \in \mathfrak{S}_A$, the  reduction of $\pi$ is the permutation in $S_n$ by replacing each letter $a_j$  by $j$.  For example, $\pi = 9\, 3\, 8\, 10\, 12\, 2\, 7$  is the permutation of the set $A=\{2,\,3,\,7,\,8,\,9,\,10,\,12\}$, its reduction is $5\, 2\, 4\, 6 \,7\, 1\, 3$, which is a permutation in  $S_7$.

 Let $\pi=\pi_1\cdots\pi_n \in S_n$. The derangement part of $\pi$, denoted $dp(\pi)$, is the reduction of the subword of non-fixed points of $\pi$. Note that $\pi_i$ is called the non-fixed points of $\pi$ if $\pi_i \neq i$. For example, let's take the following permutation in $S_9$:
\begin{equation}\label{exa-s}
\pi = 1\, 5\, 3\, 7\, 6\, 2\, 9\, 8\, 4,
\end{equation}
there are three fixed points, which are $1,3,8$ and six non-fixed points:  $5,\,7,\,6,\,2,\,9,\,4$. The reduction of non-fixed points of $\pi$ is $3\, 5\, 4\, 1\,6\, 2$, so the derangement part of $\pi$ is
\[dp(\pi)= 3\, 5\, 4\, 1\,6\, 2.\]

 Wachs \cite{Wachs-1989}   established the following relation by constructing a bijection on $S_n$  and utilizing  Garsia and Gessel's result on shuffles of permutations in \cite{Garsia-Gessel-1979}.
\begin{prop}[Wachs]\label{DSn}
Let $0\leq k\le n$ and $\sigma\in \mathcal{D}^S_k$. Then
\begin{equation}\label{DSn-eq}
    \sum_{dp(\pi)=\sigma \atop \pi \in S_n} q^{{\rm maj}(\pi)}=q^{{\rm maj}(\sigma)}{n \brack k}_q,
\end{equation}
where
\[{n \brack k}_q=\frac{[n]_q!}{[k]_q![n-k]_q!}\]
is the $q$-binomial coefficients.
\end{prop}
Chen and Xu \cite{Chen-Xu-2008}   provided an alternative bijective proof of  the  relation \eqref{DSn-eq}. To prove Theorem \ref{ThSn}, we establish the following proposition.

 \begin{prop}\label{DSn-ss}
Let $0\leq k\le n$ and $\sigma\in \mathcal{D}^S_k$, we have
\begin{equation}\label{DSn-sseq}
    \sum_{dp(\pi)=\sigma \atop \pi \in S_n} (-1)^{\ell_A(\pi)}q^{{\rm maj}(\pi)}=(-1)^{\ell_A(\sigma)}q^{{\rm maj}(\sigma)}{n \brack k}_q.
\end{equation}
\end{prop}
\pf In light of  Proposition \ref{DSn}, it suffices to show
\begin{equation} \label{snmmtt}
\ell_A(\pi) \equiv \ell_A(dp(\pi)) \pmod{2}
\end{equation}
for each permutation $\pi \in S_n$.

Let $\pi \in S_n$. Assume that there are $k$ fixed points in $\pi$, which are $i_1<\cdots<i_k$. The fixed point $i_j$  is called the $j$-th fixed point of $\pi$.  Let $\sigma=dp(\pi)$,   which clearly belongs to  $ \mathcal{D}^S_{n-k}$. Let $\pi^{(j)}$ be the reduction of the permutation obtained from $\pi^{(j-1)}$ by removing the $j$-th fixed point of $\pi$.  It is evident that $\pi^{(k)}=\sigma$, and  we assume that $\pi^{(0)}=\pi$.
We aim to demonstrate that for $1\leq j\leq k$,
\begin{equation}\label{snmm}
\ell_A(\pi^{(j-1)}) \equiv \ell_A(\pi^{(j)}) \pmod{2}.
\end{equation}
Let $\pi^{(j)}=\pi^{(j)}_1\cdots \pi^{(j)}_{n-j}$. Note that 
$i_1<\cdots<i_k$ are the fixed points of $\pi$. We next intend to  insert the $j$-th fixed point $i_j$ into $\pi^{(j)}$ to get $\pi^{(j-1)}$. Assume that $i_j-j+1=t$. Clearly, $1\leq t< n-j+2$. By definition, we see that 
$\pi^{(j-1)}$ is the permutation obtained from $\pi^{(j)}$ by 
 replacing the element $\pi^{(j)}_i$  with $\pi^{(j)}_i+1$ for   $\pi^{(j)}_i\ge t$ and inserting $t$ such that  $t$ becomes a fixed point in $\pi^{(j-1)}$. 
It is easy to check that 
\begin{equation}\label{snaa}
\ell_A(\pi^{(j-1)})-\ell_A(\pi^{(j)})=\#\{i\colon i<t, \pi^{(j)}_i\geq  t\}+\#\{i\colon i\ge t, \pi^{(j)}_i< t\}.
\end{equation}
On the other hand, it is evident that
\begin{align*}
    &\#\{i\colon i<t, \pi^{(j)}_i\geq t\}+\#\{i\colon i<t, \pi^{(j)}_i<t\}=t-1,\\[5pt]
    &\#\{i\colon i\ge t,\pi^{(j)}_i<t\}+\#\{i\colon i<t, \pi^{(j)}_i<t\}=t-1.
\end{align*}
This implies that
\begin{equation}\label{snbb}
\#\{i\colon i<t, \pi^{(j)}_i\geq  t\}=\#\{i\colon i\ge t, \pi^{(j)}_i< t\}.
\end{equation}
Combining \eqref{snaa} and \eqref{snbb}, we conclude that \eqref{snmm},  and thus,  \eqref{snmmtt} holds. Therefore, we deduce \eqref{DSn-sseq}  from \eqref{DSn-eq} and \eqref{snmmtt}. This completes the proof. \qed

With Proposition \ref{DSn-ss} at hand, we can prove Theorem \ref{ThSn} along the lines of Wachs \cite{Wachs-1989}.

\noindent{\it Proof of Theorem \ref{ThSn}:}
Summing over all derangements $\sigma \in \mathcal{D}^S_k$ and $0 \leq k \leq n$, and applying \eqref{S-1}, we can infer from Proposition \ref{DSn-ss}  that
 \begin{equation*}
    \sum_{k=0}^{n} {n \brack k}_q \sum_{\sigma \in \mathcal{D}^S_k} (-1)^{\ell_A(\sigma)} q^{{\rm maj}(\sigma)}=\left(\frac{1-q} {1+q}\right)^{\lfloor \frac{n}{2}\rfloor}[n]_q!.
\end{equation*}
Thus, Theorem \ref{ThSn} is established through the application of the $q$-binomial inversion \cite[Corollary 3.38]{Aigner-1979}.  This completes the proof. \qed

Combining Theorem \ref{dSn-Wachs} and Theorem \ref{ThSn}, we obtain the following formula for the Mahonian polynomial over the set $\mathcal{D}^{AS}_n$ of derangements in the alternating subgroup $A_n$, which immediately yields Theorem \ref{en-as} by setting $q\rightarrow 1$.
\begin{coro}\label{ThAAn} For $n\geq 1$,
    \begin{align*}
    {d}^{\,AS}_n(q)=  \sum_{\sigma \in \mathcal{D}^{AS}_n} q^{{\rm maj}(\sigma)}
    &= [n]_q!\sum^n_{k=0}\frac{(-1)^k}{[k]_q!}q^{\binom{k}{2}}\left(\frac{1}{2}+\frac{1}{2}\left(\frac{1-q} {1+q}\right)^{\lfloor \frac{n-k}{2}\rfloor}\right).\\
    \end{align*}
\end{coro}
The first six ${d}^{\,AS}_n(q)$ ($1\leq n\leq 6$) are given as follows:
    \begin{align*}
    {d}^{\,AS}_1(q)=&{d}^{\,AS}_2(q)=0;\\[5pt]
    {d}^{\,AS}_3(q)=&q+q^2;\\[5pt] {d}^{\,AS}_4(q)=&{q}^{2}+{q}^{4}+{q}^{6};\\[5pt] {d}^{\,AS}_5(q)=&q+2\,{q}^{2}+3\,{q}^{3}+4\,{q}^{4}+4\,{q}^{5}+4\,{q}^{6}+3\,{q}^{7}+2
\,{q}^{8}+{q}^{9}
;\\[5pt]
{d}^{\,AS}_6(q)=&2\,{q}^{2}+3\,{q}^{3}+8\,{q}^{4}+10\,{q}^{5}+17\,{q}^{6}+17\,{q}^{7}+
21\,{q}^{8}+17\,{q}^{9} +16\,{q}^{10}+9\,{q}^{11}\\[5pt] &\quad+7\,{q}^{12}+2\,{q}^{
13}+{q}^{14}.
    \end{align*}

\section{The signed Mahonian polynomials over $\mathcal{D}^B_n$}

This section is focused on exploring the signed Mahonian polynomials within the context of signed derangements in the hyperoctahedral group $B_n$. Notably, the hyperoctahedral group $B_n$ is recognized as a Coxeter group of type $B$  with the set of generators $S = \{s^B_0, s_1 ,\ldots,s_{n-1} \}$, where 
$s^B_0=[-1,2,3,\ldots,n]$ and  $s_i=[1,2,\ldots,i-1,i+1,i,i+2,\ldots,n]$ \text{for} $1\leq i \leq n-1,$
see \cite[Proposition 8.1.3]{Bjorner-Brenti-2005}. Any element $\sigma$ of $B_n$ can be represented by a signed word $\sigma_1 \sigma_2\cdots \sigma_n$ of length $n$, where $\sigma _i \in \{-n,\ldots, -1, 1,\ldots, n\}$ and $|\sigma_1| \cdots |\sigma_n|$ is a permutation of the set $[n]$.   It is easy to see that $|B_n|=2^n n!$.

Similarly,   one can employ direct combinatorial approaches to calculate the length $\ell_B(\sigma)$ and the descent set ${\rm Des}_B(\sigma)$ (see \cite[Proposition 8.1.1 and Proposition 8.1.2]{Bjorner-Brenti-2005} and \cite[Proposition 3.1 and Corollary 3.2]{Brenti-1994}).

For $\sigma=\sigma_1 \cdots \sigma_n\in B_n$, let
\begin{align*}
&{\rm Neg}(\sigma):=\{i\in[n]: \sigma_i<0 \},\\[5pt]
&{\rm neg}(\sigma):=\#{\rm Neg}(\sigma).
\end{align*}
The length $\ell_B(\sigma)$ and the descent set ${\rm Des}_B(\sigma)$ are given by
\begin{equation}\label{def-lb}
 \ell_B(\sigma)={\rm inv}(\sigma_1 \cdots \sigma_n)-\sum_{i \in {\rm Neg}(\sigma)}\sigma_i.
\end{equation}
and
\begin{equation*}
 {\rm Des}_B(\sigma)=\{0\leq i \leq n-1 \colon  \sigma_i >\sigma_{i+1}\},
\end{equation*}
where we assume that $\sigma_0=0$.  Here we adopt the following order of the elements of $B_n$ to compute ${\rm inv}(\sigma_1 \cdots \sigma_n)$ as define in \eqref{defi-inv}.
\begin{equation}\label{naord}
    -n< -(n-1)< \cdots < -1 < 0 < 1 < 2 < \cdots< n.
\end{equation}
For example,
 let  $\sigma=\bar{3}\,1\,\bar{6}\,2\,\bar{4}\,\bar{5}\in B_6$, we see that
\[
\ell_B(\sigma)={\rm inv}(\sigma)+3+6+4+5=27, \quad \text{and} \quad  {\rm des}_B(\sigma)=4.
\]

The first major index ${\rm maj}_B$ on the hyperochaherdral group $B_n$  was introduced by
 Adin and Roichman \cite[Theorem 3]{Adin-Roichman-2001}, referred to as the  flag major index, denoted ${\rm fmaj}$:
\begin{equation}\label{defi-fmaj}
    {\rm fmaj}(\sigma )=2{\rm maj}(\sigma_1 \cdots \sigma_n)+{\rm neg}(\sigma),
\end{equation}
where ${\rm maj}(\sigma_1 \cdots \sigma_n)$, as defined in \eqref{def-maj},  is computed  using the following order on $\mathbf{Z}$:
\begin{equation} \label{neword}
    -1\prec -2\prec \cdots \prec -n  \prec 0 \prec 1 \prec 2 \prec \cdots\prec n
\end{equation}
instead of the usual ordering \eqref{naord}.

For example, let $\sigma=\bar{3}\,1\,\bar{6}\,2\,\bar{4}\,\bar{5}\in B_6$.  The corresponding flag major index is given by
\[{\rm fmaj}(\sigma)=(2+4)\times 2+4=16.\]

Adin and Roichman \cite{Adin-Roichman-2001} proved that the flag major index is a  Mahonian statistic  of the hyperoctahedral group $B_n$. Specifically, they established the following result: 
\begin{thm}[Adin-Roichman]\label{BM-1} For $n\geq 1$,
\begin{equation*}
 \sum_{\sigma \in B_n}q^{{\rm fmaj}(\sigma)}=\sum_{\sigma \in B_n}q^{\ell_B(\sigma)}=[2]_{q}[4]_{q}\cdots[2n]_{q}.
\end{equation*}
\end{thm}
Adin, Gessel, and Roichman \cite{Adin-Gessel-Roichman-2005} derived the following formula for the signed Mahonian polynomial over the hyperochaherdral group $B_n$:
\begin{thm}[Adin-Gessel-Roichman]\label{BMS-1} For $n\geq 1$,
\begin{eqnarray*}
 \sum_{\sigma \in B_n}(-1)^{\ell_B(\sigma)}q^{{\rm fmaj}(\sigma)}&=&[2]_{-q}[4]_{q}\cdots[2n]_{(-1)^n q}\\[5pt]
 &=&\left(\frac{1-q}{1+q}\right)^{\lfloor \frac{n-1}{2}\rfloor}[2]_{q}[4]_{q}\cdots[2n]_{q}.
\end{eqnarray*}
\end{thm}
The following formula for the Mahonian polynomial over $\mathcal{D}^B_n$ is due to  Chow \cite{Chow-2006}.
\begin{thm}[Chow]\label{dBn} For $n\geq 1$,
\begin{equation}
  d^{B}_n(q)=\sum_{\sigma \in \mathcal{D}^B_n} q^{{\rm fmaj}(\sigma)}=[2]_q[4]_q\cdots[2n]_q\sum^n_{k=0}\frac{(-1)^kq^{2\binom{k}{2}}}{[2]_q[4]_q\cdots[2k]_q}.
\end{equation}
\end{thm}
The second result of this paper is  the following formula for   the signed Mahonian polynomial over $\mathcal{D}^B_n$:
\begin{thm}\label{ThBn} For $n\geq 1$,
    \begin{equation*}
    \overline{d}^{B}_n(q)=\sum_{\sigma \in \mathcal{D}^B_n} (-1)^{\ell_B(\sigma)} q^{{\rm fmaj}(\sigma)}=[2]_q[4]_q\cdots[2n]_q\sum^n_{k=0}\frac{(-1)^kq^{2\binom{k}{2}}}{[2]_q[4]_q\cdots[2k]_q}\left(\frac{1-q} {1+q}\right)^{\lfloor \frac{n-k-1}{2}\rfloor}.
\end{equation*}
\end{thm}

Similarly, let us first recall the definition of the derangement part of the signed permutation $\pi \in B_n$.
 Let $A=\{0<a_1<a_2<\cdots <a_n\}$   and let
 $\pi=\pi_1\cdots \pi_n$ be a signed permutation so that  $|\pi_1||\pi_2|\cdots |\pi_n|$ is a permutation in $\mathfrak{S}_A$.  The reduction of $\pi$ is the signed permutation in $B_n$ by replacing each letter $\pi_i=a_j$  by $({\rm sgn}~\pi_i)j$.   For example, $\pi = 2\, 5\, \bar{3}\, 8\, \bar{9}$  is the signed permutation of the set $A=\{2,\,3,\,5,\,8,\,9\}$, its reduction is $1\, 3\, \bar{2}\, 4 \, \bar{5}$, which is a signed permutation in  $B_5$.

  The derangement part of $\pi$, denoted $dp(\pi)$, is the reduction of the signed subword of non-fixed points of $\pi$. Recall that $\pi_i$ is called the non-fixed points of $\pi=\pi_1\cdots \pi_n$ if $\pi_i \neq i$. For example, consider the permutation
\[\pi= 1\, 6\, \bar{3}\, 5\, 8\, 2\, 7\, \bar{4} \in B_8.\]
It has two fixed points: $1,\,7$ and six non-fixed points: $6,  \bar{3}, 5, 8, 2, \bar{4}$. The reduction of non-fixed points of $\pi$ yields $5\,  \bar{2}\, 4\, 6\, 1\, \bar{3}$. Consequently,  the derangement part of $\pi$ is
\[dp(\pi)= 5\,  \bar{2}\, 4\, 6\, 1\, \bar{3} \in B_6.\]
By extending Wachs' approach to the context of signed permutations in $B_n$,  Chow \cite{Chow-2006} established the following relation, which is the main ingredient in the derivation of Theorem  \ref{dBn}.
\begin{prop}[Chow]\label{DBn}
Let  $0\le k \le n$ and $\sigma \in {\mathcal D}^B_k$, we have
\begin{equation}\label{DBn-eq}
    \sum_{dp(\pi)=\sigma  \atop \pi \in B_n} q^{{\rm fmaj}(\pi)}=q^{{\rm fmaj}(\sigma )}{n \brack k}_{q^2}.
\end{equation}
\end{prop}
In a similar vein,  to prove Theorem \ref{ThBn}, we need to establish the following proposition.
\begin{prop}\label{DBsn}
Let $0\le k \le n$ and $\sigma \in \mathcal{D}^B_k$, we have
\begin{equation}\label{DBn-sseq}
    \sum_{dp(\pi)=\sigma \atop \pi \in B_n} (-1)^{\ell_B(\pi)} q^{{\rm fmaj}(\pi)}=(-1)^{\ell_B(\sigma)}q^{{\rm fmaj}(\sigma)}{n \brack k}_{q^2}.
\end{equation}
\end{prop}
\pf Building upon Proposition \ref{DBn}, it suffices to show that
\begin{equation} \label{bnmmtt}
\ell_B(\pi) \equiv \ell_B(dp(\pi)) \pmod{2}
\end{equation}
for each permutation $\pi \in B_n$.

Let $\pi \in B_n$ and suppose that there are $k$ fixed points in $\pi$, which are $i_1<\cdots<i_k$. The fixed point $i_j$  is called the $j$-th fixed point of $\pi$.  Let $\sigma=dp(\pi)$. It is clear that $\sigma \in \mathcal{D}^B_{n-k}$. Let $\pi^{(j)}$ be the reduction of the permutation obtained from $\pi^{(j-1)}$ by removing the $j$-th fixed point of $\pi$. Here we assume that $\pi^{(0)}=\pi$. Clearly, $\pi^{(k)}=\sigma$.
We proceed to show that for $1\leq j\leq k$,
\begin{equation}\label{bnmm}
\ell_B(\pi^{(j-1)}) \equiv \ell_B(\pi^{(j)}) \pmod{2}.
\end{equation}
Let $\pi^{(j)}=\pi^{(j)}_1\cdots \pi^{(j)}_{n-j}$. Recall  that 
$i_1<\cdots<i_k$ are the fixed points of $\pi$. We aim to insert the $j$-th fixed point $i_j$ into $\pi^{(j)}$ to get $\pi^{(j-1)}$. More precisely, assume that $i_j-j+1=t$, where clearly $1\leq t< n-j+2$. By definition, we see that 
$\pi^{(j-1)}$ is the signed permutation obtained from $\pi^{(j)}$ by 
 replacing the element $\pi^{(j)}_i$  with ${\rm sgn }\
  \pi^{(j)}_i(|\pi^{(j)}_i|+1)$ for   $|\pi^{(j)}_i|\ge t$ and inserting $t$ such that  $t$ becomes a fixed point in $\pi^{(j-1)}$. 
To prove \eqref{bnmm}, we consider the following eight sets:
\begin{align*}
   M^{> ,+}_{<t}&=\{i\colon i<t,|\pi^{(j-1)}_i|> t, \pi^{(j-1)}_i>0 \};\\[5pt]
   M^{>,-}_{<t}&=\{i\colon i<t,|\pi^{(j-1)}_i|> t, \pi^{(j-1)}_i<0 \};\\[5pt]
   M^{<,+}_{<t}&=\{i\colon i<t,|\pi^{(j-1)}_i|< t, \pi^{(j-1)}_i>0 \};\\[5pt]
   M^{<,-}_{<t}&=\{i\colon i<t,|\pi^{(j-1)}_i|< t, \pi^{(j-1)}_i<0 \};\\[5pt]
   M^{>,+}_{>t}&=\{i\colon i>t,|\pi^{(j-1)}_i|> t, \pi^{(j-1)}_i>0 \};\\[5pt]
   M^{>,-}_{>t}&=\{i\colon i>t,|\pi^{(j-1)}_i|> t, \pi^{(j-1)}_i<0 \};\\[5pt]
   M^{<,+}_{>t}&=\{i\colon i>t,|\pi^{(j-1)}_i|< t, \pi^{(j-1)}_i>0 \};\\[5pt]
   M^{<,-}_{>t}&=\{i\colon i>t,|\pi^{(j-1)}_i|< t, \pi^{(j-1)}_i<0 \}.
\end{align*}
By definition, it is not difficult to show that
\begin{align*}
   {\rm inv}(\pi^{(j-1)})-{\rm inv}(\pi^{(j)})&= \#M^{> ,+}_{<t}+\# M^{>,-}_{>t}+\# M^{<,+}_{>t}+\#M^{<,-}_{>t}
\end{align*}
and
\begin{align*}
   -\sum_{i \in {\rm Neg}(\pi^{(j-1)})}\pi^{(j-1)}_i+\sum_{i \in {\rm Neg}(\pi^{(j)})}\pi^{(j)}_i&=\# M^{>, -}_{<t}+\# M^{>,-}_{>t}.
\end{align*}
Hence, we derive from \eqref{def-lb} that
\begin{align}\label{bnaa}
   \ell_B(\pi^{(j-1)}) - \ell_B(\pi^{(j)})&={\rm inv}(\pi^{(j-1)})-\sum_{i \in {\rm Neg}(\pi^{(j-1)})}\pi^{(j-1)}_i-{\rm inv}(\pi^{(j)})+\sum_{i \in {\rm Neg}(\pi^{(j)})}\pi^{(j)}_i \nonumber \\[5pt]
&=\# M^{> ,+}_{<t}+ \# M^{>,-}_{<t}+2 \# M^{>,-}_{>t}+ \# M^{<,+}_{>t} +\# M^{<,-}_{>t}.
\end{align}
On the other hand, it is easy to check that
\begin{align*}    \# M^{>,+}_{<t}+\# M^{>,-}_{<t}+\# M^{<,+}_{<t}+\# M^{<,-}_{<t}&=\# \{i:i<t\}=t-1,\\[5pt]    \# M^{<,+}_{>t}+\# M^{<,-}_{>t}+\# M^{<,+}_{<t}+\# M^{<,-}_{<t}&=\# \{i:|\pi^{(j-1)}_i|<t\}=t-1.
\end{align*}
Consequently, 
\begin{equation}    \label{bnbb}
\# M^{>,+}_{<t}+\# M^{>,-}_{<t}=\# M^{<,+}_{>t}+\# M^{<,-}_{>t}.
\end{equation}
By substituting \eqref{bnbb} into \eqref{bnaa}, we obtain \eqref{bnmm}, establishing the validity of \eqref{bnmmtt}. Consequently, we demonstrate Proposition \ref{DBsn} by combining \eqref{DBn-eq} and \eqref{bnmmtt}. This completes the proof. \qed

We  are ready to  prove Theorem \ref{ThBn} with the aid of  Proposition \ref{DBsn}.

\noindent{\it Proof of Theorem \ref{ThBn}:}
Summing \eqref{DBn-sseq} in Proposition  \ref{DBsn} over all derangements $\sigma \in \mathcal{D}^B_k$ for  $0 \leq k \leq n$ and  applying Theorem \ref{BMS-1},  we obtain  that
 \begin{equation*}
    \sum_{k=0}^{n} {n \brack k}_{q^2} \sum_{\sigma \in \mathcal{D}^B_k} (-1)^{\ell_B(\sigma)} q^{{\rm fmaj}(\sigma)}=\left(\frac{1-q}{1+q}\right)^{\lfloor \frac{n+1}{2}\rfloor}[2]_{q}[4]_{q}\cdots[2n]_{q},
\end{equation*}
which leads to Theorem \ref{ThBn} by utilizing the $q$-binomial inversion \cite[Corollary 3.38]{Aigner-1979}.  This completes the proof. \qed

Combining Theorem \ref{dBn} and  Theorem \ref{ThBn}, we derive the following formula for the Mahonian polynomial over the set of derangements of type $B$  with even length. This formula immediately leads to   Theorem \ref{en-ab} when substituting $q\rightarrow 1$.

\begin{coro}\label{ThABn} For $n\geq 1$,
    \begin{align*}
    {d}^{\,AB}_n(q)=  \sum_{\sigma \in \mathcal{D}^{AB}_n} q^{{\rm fmaj}(\sigma)}    &= {[2]_q[4]_q\cdots[2n]_q}\sum^n_{k=0}\frac{(-1)^kq^{2\binom{k}{2}}}{[2]_q[4]_q\cdots[2k]_q}\left(\frac{1}{2}+\frac{1}{2}\left(\frac{1-q} {1+q}\right)^{\lfloor \frac{n-k+1}{2}\rfloor}\right).\\
    \end{align*}
\end{coro}
Below are the first five terms of ${d}^{\,AB}_n(q)$ for $1\leq n\leq 5$:
    \begin{align*}
    {d}^{\,AB}_1(q)=&0;\\[5pt]
       {d}^{\,AB}_2(q)=&q+q^2+q^3;\\[5pt]
       {d}^{\,AB}_3(q)=&2\,{q}^{2}+{q}^{3}+3\,{q}^{4}+2\,{q}^{5}+2\,{q}^{6}+2\,{q}^{7}+{q}^{8}
+{q}^{9}
;\\[5pt]
        {d}^{\,AB}_4(q)=&q+2\,{q}^{2}+5\,{q}^{3}+6\,{q}^{4}+10\,{q}^{5}+10\,{q}^{6}+14\,{q}^{7}
+13\,{q}^{8}+14\,{q}^{9}+12\,{q}^{10}+10\,{q}^{11}
\\[5pt]
        & +9\,{q}^{12}+5\,{q}^
{13}+4\,{q}^{14}+{q}^{15}+{q}^{16};\\[5pt]
       {d}^{\,AB}_5(q)=&3\,{q}^{2}+5\,{q}^{3}+14\,{q}^{4}+20\,{q}^{5}+34\,{q}^{6}+44\,{q}^{7}+
61\,{q}^{8}+73\,{q}^{9}+87\,{q}^{10}+97\,{q}^{11}+103\,{q}^{12}
\\[5pt]
        &+106\,{
q}^{13}+101\,{q}^{14}+96\,{q}^{15}+83\,{q}^{16}+72\,{q}^{17}+56\,{q}^{
18}+43\,{q}^{19}+29\,{q}^{20}+19\,{q}^{21} \\[5pt]
        &+11\,{q}^{22}+5\,{q}^{23}+2
\,{q}^{24}.
    \end{align*}

\section{The (signed) Mahonian polynomials over $\mathcal{D}^D_n$}

 This section aims to study the Mahonian polynomials  and the signed Mahonian polynomials over the set of derangements in the   even-signed permutation group $D_n$. The even-signed permutation group $D_n$ is a Coxeter group of type $D$ with the set of generators $S = \{s^D_0, s_1 ,\ldots,s_{n-1} \}$, where $s^D_0=[-2,-1,3,\ldots,n]$ and $s_i=[1,2,\ldots,i-1,i+1,i,i+2,\ldots,n]$ for $1\leq i \leq n-1$  (see \cite[Proposition 8.2.3]{Bjorner-Brenti-2005}).

Let $\sigma=\sigma_1\cdots \sigma_n \in D_n$,  we see that $|\sigma_1|\cdots |\sigma_n|$ is a permutation of $[n]$ and there are even number of negative entries in $\sigma$. Likewise,  the length $\ell_D(\sigma)$ and the descent set ${\rm Des}_D(\sigma)$ can be described in a combinatorial way. For the detailed proofs, please refer to   \cite[Proposition 8.2.1 and Proposition 8.2.2]{Bjorner-Brenti-2005}.

The length $\ell_D(\sigma)$ can be computed as follows:
\begin{align*}
 \ell_D(\sigma)={\rm inv}(\sigma_1\cdots \sigma_n)-\sum_{i \in {\rm Neg}(\sigma)}(\sigma_i+1)=\ell_B(\sigma)- {\rm neg}(\sigma),
\end{align*}
and the descent set ${\rm Des}_D(\sigma)$ is given by
\begin{equation*}
 {\rm Des}_D(\sigma)=\{0\leq i \leq n-1 \colon  \sigma_i >\sigma_{i+1}\},
\end{equation*}
where  $\sigma_0=-\sigma_2$. We adhere to the  order of the elements of $B_n$ given in \eqref{naord} to compute ${\rm inv}(\sigma_1\cdots \sigma_n)$ as defined in \eqref{defi-inv}.
For example, let  $\sigma=\bar{3}\,1\,\bar{6}\,2\,\bar{4}\,\bar{5}\in D_6$, we see that
${\rm inv}(\sigma)=9$,
\[
\ell_D(\sigma)=23, \quad {\rm and} \quad {\rm des}_D(\sigma)=4.
\]
Biagioli and Caselli \cite{Biagioli-Caselli-2004} introduced the D-major index, denoted as ${\rm Dmaj}$, which has been proven to be a major index  ${\rm maj}_D$ on the even-signed permutation group $D_n$.

For any $\sigma=\sigma_1\cdots \sigma_n \in  D_n $ let
\begin{equation*}
   |\sigma|_n:=\sigma_1\cdots \sigma_{n-1}|\sigma_n|.
\end{equation*}
The  D-major index of $\sigma$ is defined as
\begin{equation*}
    {\rm Dmaj}(\sigma):={\rm fmaj}(|\sigma|_n).
\end{equation*}

For $\sigma=\bar{3}\,1\,\bar{6}\,2\,\bar{4}\,\bar{5}\in D_6$, we see that $|\sigma|_n=\bar{3}\,1\,\bar{6}\,2\,\bar{4}\, {5}$, and so
\[ {\rm Dmaj}(\sigma):={\rm fmaj}(|\sigma|_n)=2{\rm maj}(\bar{3}\,1\,\bar{6}\,2\,\bar{4}\, {5})+3=15.\]
It should be stressed that ${\rm maj}(\bar{3}\,1\,\bar{6}\,2\,\bar{4}\, {5})$ is computed using the order given by \eqref{neword}.

The following result due to Biagioli and Caselli \cite{Biagioli-Caselli-2004} demonstrates  that  ${\rm Dmaj}(\sigma)$ is a Mahonian statistic on  the   even-signed permutation group $D_n$.
\begin{thm}[Biagioli-Caselli] \label{bc-mahonian} For $n\geq 1$,
\begin{equation}
\sum _ { \sigma \in D _ { n } } q ^ { \ell_D ( \sigma ) }  = \sum _ { \sigma \in D _ { n } } q ^ { \operatorname{Dmaj} ( \sigma ) }=[2]_q[4]_q\cdots [2n-2]_q [n]_q.
\end{equation}
\end{thm}
Biagioli  \cite{Biagioli-2006} obtained the following formula for the signed Mahonian polynomial over the   even-signed permutation group $D_n$ equipped with the D-major index.
\begin{thm}[Biagioli] \label{sigMahothm} For $n\geq 1$,
\begin{equation}\label{sigMaho}
 \sum_{\sigma\in D_n}(-1)^{\ell_D(\sigma)}q^{{\rm Dmaj} ( \sigma )}=[2]_{-q}[4]_q\cdots[2n-2]_{(-1)^{n-1}q}[n]_q.
\end{equation}
\end{thm}
Recently, Chow \cite{Chow-2023} pioneered the exploration of the $q$-derangement polynomial in $D_n$ using the major index ${\rm maj}_A(\sigma)$. With the aid of a result due to Foata and Han \cite{Foata-Han-2008}, he obtained the following formula:
 \begin{thm}[Chow]\label{ThDnchow}  For $n\geq 1$,
$$\sum _ { \sigma \in \mathcal{D}^D_n }  q ^ { {\rm maj}_A(\sigma) }  = \sum _ { k = 0 } ^ { n } ( - 1 ) ^ { k } q ^ { \tbinom{k}{2} } 2 ^ { n - 1 - k } \frac { [ n ] _ { q } ! } { [ k ] _ { q } ! } + \frac { ( - 1 ) ^ { n } q ^ {\tbinom{n}{2} } } { 2 },$$
where ${\rm maj}_A(\sigma)={\rm maj}(\sigma_1\sigma_2\cdots \sigma_n)$ as defined in \eqref{def-maj}  is computed according to the  order given in \eqref{naord}.
 \end{thm}

The third result of this paper is to establish the following formula for the  Mahonian polynomials over the set of derangements in $D_n$  equipped with the D-major index.
\begin{thm}\label{ThDn+} For $n\geq 2$,
    \begin{align}\label{ThDn+-eq}
     d^{D}_n(q)=\sum_{\sigma \in \mathcal{D}^{D}_n}q^{\operatorname { Dmaj }(\sigma)}&=[2]_q[4]_q\cdots[2n-2]_q\sum^{n-2}_{k=0}\frac{(-1)^kq^{2\binom{k}{2}}}{[2]_q[4]_q\cdots[2k]_q} \notag\\
    &\quad \times \left(q^{2k+1}[n-1-k]_q+ \frac{1}{2}\left(1-\left(\frac{1-q} {1+ q}\right)^{ n-k-1}\right) \right).
\end{align}
\end{thm}
Below are the first five terms of $d^{D}_n(q)$ for $1\leq n\leq 5$:
    \begin{align*}
    {d}^{\,D}_1(q)=&0;\\[5pt]
       {d}^{\,D}_2(q)=&2q+q^2;\\[5pt]
       {d}^{\,D}_3(q)=&2\,q+3\,{q}^{2}+4\,{q}^{3}+3\,{q}^{4}+2\,{q}^{5};\\[5pt]
        {d}^{\,D}_4(q)=&2\,q+4\,{q}^{2}+12\,{q}^{3}+12\,{q}^{4}+20\,{q}^{5}+16\,{q}^{6}+20\,{q
}^{7}\\[5pt]
&+13\,{q}^{8}+10\,{q}^{9}+5\,{q}^{10}+2\,{q}^{11}+{q}^{12} ;\\[5pt]
{d}^{\,D}_5(q)=&2\,q+5\,{q}^{2}+21\,{q}^{3}+26\,{q}^{4}+61\,{q}^{5}+62\,{q}^{6}+108\,{
q}^{7}+100\,{q}^{8}
 \end{align*}
 \begin{align*}
          &\quad \quad +141\,{q}^{9}+118\,{q}^{10}+136\,{q}^{11}+105\,{q}^
{12}+99\,{q}^{13}+69\,{q}^{14}+52\,{q}^{15}\\[5pt]
        &\quad \quad+31\,{q}^{16}+17\,{q}^{17}+
8\,{q}^{18}+3\,{q}^{19}.
    \end{align*}
Like the symmetric group $S_n$ and the   hyperoctahedral group $B_n$,  we also establish the formula for the  signed  Mahonian polynomials over the set of derangements in $D_n$.
\begin{thm}\label{ThDn-} For $n\geq 2$,
    \begin{align}\label{ThDn--eq}
   \overline{d}^{D}_n(q)&= \sum_{\sigma \in \mathcal{D}^{D}_n}(-1)^{\ell_D(\sigma)}q^{\operatorname { Dmaj }(\sigma)} \nonumber \\[3pt]
   &=[2]_q[4]_q\cdots[2n-2]_q\sum^{n-2}_{k=0}\frac{(-1)^{n-1}q^{2\binom{k}{2}}}{[2]_q[4]_q\cdots[2k]_q} \left(\frac{1-q}{1+ q}\right)^{\lfloor \frac{n-k}{2}\rfloor}  \nonumber\\[5pt]
     &\quad \times\frac{2q^{2k+1}(1+(-1)^{n-k}q^{n-k-1})-(1+(-1)^{n-k})q}{2(1-q)}.
\end{align}
\end{thm}
In a similar vein,  the combination of Theorem \ref{ThDn+} and Theorem \ref{ThDn-} results in the formula for the Mahonian polynomial over the set of derangements of type $D$ with an even length.  This formula directly implies Theorem \ref{en-ad} by setting $q\rightarrow 1$.
\begin{coro}\label{ThADn} For $n\geq 2$,
\begin{align*}
   d^{AD}_n(q)&=\sum_{\sigma \in \mathcal{D}^{AD}_n} q^{{\rm Dmaj}(\sigma)} \\
   &=  [2]_q[4]_q\cdots[2n-2]_q\sum^{n-2}_{k=0}\frac{(-1)^kq^{k^2-k}}{[2]_q[4]_q\cdots[2k]_q} \\[5pt]
     &\quad \times  \left(\frac{1}{2}q^{2k+1}[n-1-k]_q+ \frac{1}{4}\left(1-\left(\frac{1-q} {1+ q}\right)^{ n-k-1}\right)\right.\\[5pt]
     &\quad \left.+\left(\frac{1-q}{1+ q}\right)^{\lfloor \frac{n-k}{2}\rfloor}\left(\frac{2q^{2k+1}(1+(-1)^{n-k}q^{n-k-1})-(1+(-1)^{n-k})q}{4(1-q)}\right)\right).
\end{align*}
\end{coro}
Below are the first five terms of ${d}^{\,AD}_n(q)$ for $1\leq n\leq 5$:
    \begin{align*}
    {d}^{\,AD}_1(q)=&0;\\[5pt]
       {d}^{\,AD}_2(q)=&q;\\[5pt]
       {d}^{\,AD}_3(q)=&q+2\,{q}^{2}+2\,{q}^{3}+2\,{q}^{4}+{q}^{5};\\[5pt]
        {d}^{\,AD}_4(q)=&q+2\,{q}^{2}+6\,{q}^{3}+5\,{q}^{4}+10\,{q}^{5}+7\,{q}^{6}+10\,{q}^{7}+
6\,{q}^{8}+5\,{q}^{9}+3\,{q}^{10}+{q}^{11}+{q}^{12}
;\\[5pt]
       {d}^{\,AD}_5(q)=&q+3\,{q}^{2}+11\,{q}^{3}+14\,{q}^{4}+31\,{q}^{5}+32\,{q}^{6}+54\,{q}^{
7}+51\,{q}^{8}+70\,{q}^{9}+59\,{q}^{10}\\[5pt]
&+67\,{q}^{11}+52\,{q}^{12}+49\,
{q}^{13}+34\,{q}^{14}+26\,{q}^{15}+15\,{q}^{16}+9\,{q}^{17}+4\,{q}^{18
}+2\,{q}^{19}.
    \end{align*}
The remainder of this paper is dedicated to  the proofs of Theorem \ref{ThDn+} and Theorem \ref{ThDn-}.
To  this end,   let us recall the strategy of Biagioli and Caselli  \cite{Biagioli-Caselli-2004} in the proof of Theorem \ref{bc-mahonian}.  Let \[\Delta_ { n }=\{\gamma=\gamma_1\cdots \gamma_n \in B _ { n }, \gamma_n>0\}.\]
 There is a bijection $\phi \colon D_n \rightarrow \Delta_n$ such that for all $\gamma=\gamma_1 \cdots \gamma_n\in  D_n$, we have
 \[\phi(\gamma)=\gamma_1 \gamma_2 \cdots |\gamma_n|.\]
It is evident that
 \[{\rm Dmaj}(\gamma)={\rm fmaj}(\phi(\gamma)).\]
Consequently,  the proof of Theorem \ref{bc-mahonian} is equivalent to demonstrating that
\begin{equation}\label{delta}
\sum _ { \gamma \in \Delta_ { n } } q ^ { \operatorname{\rm fmaj} ( \gamma ) }=\sum _ { \gamma \in D_ { n } } q ^ { \operatorname{\rm Dmaj} (\gamma) }=[2]_q[4]_q\cdots [2n-2]_q [n]_q.
\end{equation}
In fact, as Biagioli pointed out in \cite{Biagioli-2006}, the results concerning the D-maj index on the set $D_n$ are typically translated to investigate the results related to the fmaj index in the set $\Delta_n$.  This approach  will make definitions and arguments more natural and transparent. For example,   the key component in deriving the formula \eqref{sigMaho}  for the signed Mahonian polynomial over $D_n$ is   the following formula  for the signed Mahonian polynomial over $\Delta_n$ (see  \cite[Corollary 4.3 and Corollary 4.4]{Biagioli-2006} for the proof).
\begin{prop}[Biagioli] For $n\geq 1$,
\begin{equation}\label{deltas}
\sum _ { \gamma \in \Delta _ { n } } (-1)^{\ell_B(\gamma)} q ^ { \operatorname { fmaj } ( \gamma ) }=[2]_{-q}[4]_q\cdots[2n-2]_{(-1)^{n-1}q}[n]_{(-1)^nq}.
\end{equation}
\end{prop}
Based on this reason, it becomes essential for us to consider the set  $\mathcal{D}^{\Delta}_n$  of signed derangements in $\Delta_n$. As an illustration,
\[\mathcal{D}^{\Delta}_2=\{2\,1,\bar{2}\,1\}.\]
It is unsurprising that the formula for the (signed) Mahonian polynomials over $\mathcal{D}^{\Delta}_n$ becomes crucial  in the proofs of  Theorem \ref{ThDn+} and Theorem \ref{ThDn-}. More precisely, we have the following relation:

\begin{lem}  \label{relBD-thm} Let $\varepsilon=\pm1$ and for $n\geq 2$,
\begin{align}  \label{relBD}
 \sum_{\sigma\in \mathcal{D}^{D}_n}\varepsilon^{\ell_D(\sigma)}q^{\operatorname { Dmaj }(\sigma)}&=\sum_{\sigma\in \mathcal{D}^{\Delta}_n}\varepsilon^ { \ell _ { B } ( \sigma ) } (\varepsilon q)^{\operatorname { fmaj }(\sigma)} \nonumber\\[5pt]
 &+\frac{\varepsilon}{2}\left(\sum_{\sigma\in \mathcal{D}^B_{n-1} }\varepsilon^{\ell_B(\sigma)}q^{\operatorname { fmaj }(\sigma)}-\sum_{\sigma \in \mathcal{D}^B_{n-1} }\varepsilon^{\ell_B(\sigma)}(-q)^{\operatorname { fmaj }(\sigma)}\right).
 \end{align}
\end{lem}

\pf By definition, we have
\begin{equation} \label{relBD-a}
    \sum_{\sigma\in \mathcal{D}^{D}_n}\varepsilon^{\ell_D(\sigma)}q^{\operatorname { Dmaj }(\sigma)}=\sum_{\sigma\in \mathcal{D}^{D}_n \atop \sigma_n\neq -n}\varepsilon^{\ell_D(\sigma)}q^{\operatorname { Dmaj }(\sigma)}+\sum_{\sigma\in \mathcal{D}^D_n \atop \sigma_n=-n}\varepsilon^{\ell_D(\sigma)}q^{\operatorname { Dmaj }(\sigma)}.
\end{equation}
Let $\sigma=\sigma_1\cdots \sigma_n\in \mathcal{D}^D_n$ with $ \sigma_n=-n$, define
\[\sigma^{\prime}=\sigma_1\cdots\sigma_{n-1},\]
which is evidently a signed derangement in $\mathcal{D}^B_{n-1}$ with an odd number of negative entries. Additionally, we have
\[\ell_D(\sigma)\equiv \ell_B(\sigma)=\ell_B(\sigma^{\prime})+(2n-1) \pmod{2}  \quad \text{and} \quad  {\rm Dmaj}(\sigma)={\rm fmaj}(\sigma^{\prime}).\]
Moreover, this process is reversible. Therefore, we derive that
\begin{align}\label{relBD-b}
\sum_{\sigma\in \mathcal{D}^D_n \atop \sigma_n=-n}\varepsilon^{\ell_D(\sigma)}q^{\operatorname { Dmaj }(\sigma)}&=\sum_{\sigma'\in \mathcal{D}^B_{n-1} \atop {\rm neg}(\sigma')\equiv 1 \pmod{2}}\varepsilon^{\ell_B(\sigma')+1}q^{\operatorname { fmaj }(\sigma')}\nonumber \\[5pt]
&=\frac{\varepsilon}{2}\left(\sum_{\sigma'\in \mathcal{D}^B_{n-1} }\varepsilon^{\ell_B(\sigma')}q^{\operatorname { fmaj }(\sigma')}-\sum_{\sigma'\in \mathcal{D}^B_{n-1} }\varepsilon^{\ell_B(\sigma')}(-q)^{\operatorname { fmaj }(\sigma')}\right).
\end{align}
 Let $\sigma=\sigma_1\cdots \sigma_n\in \mathcal{D}^D_n$ with $ \sigma_n\neq -n$, define
\[|\sigma|_n=\sigma_1\cdots\sigma_{n-1}|\sigma_{n}|.\]
It is evident that $|\sigma|_n \in \mathcal{D}^{\Delta}_n$ such that
${\rm Dmaj}(\sigma)={\rm fmaj}(|\sigma|_n)$.  Moreover, it is easy to check that  $\ell_D(\sigma)\equiv \ell_B(|\sigma|_n)+\chi(\sigma_n<0) \pmod{2}$ and ${\rm fmaj}(|\sigma|_n)\equiv \chi(\sigma_n<0) \pmod{2}$, where $\chi(\sigma_n<0)=1$ if $\sigma_n<0$, and $0$ otherwise. Hence we conclude that
$$\ell_D(\sigma)\equiv \ell_B(|\sigma|_n)+ {\rm fmaj}(|\sigma|_n) \pmod{2}.$$
Moreover, this process is reversible. Consequently, we deduce that
\begin{align}\label{relBD-c}
\sum_{\sigma\in \mathcal{D}^D_n \atop \sigma_n\ne -n}\varepsilon^{\ell_D(\sigma)}q^{\operatorname { Dmaj }(\sigma)}&=\sum_{|\sigma|_n\in \mathcal{D}^\Delta_{n}}\varepsilon^{\ell_B(|\sigma|_n)}{(\varepsilon q)}^{\operatorname { fmaj }(|\sigma|_n)}.
\end{align}
Applying \eqref{relBD-b} and \eqref{relBD-c} to \eqref{relBD-a}, we get  \eqref{relBD}. This completes the proof. \qed

In light of Lemma \ref{relBD-thm}, along with the application of Theorem \ref{dBn} and Theorem \ref{ThBn}, it becomes apparent that the demonstrations of both Theorem \ref{ThDn+} and Theorem \ref{ThDn-} hinge on deducing the formula for the (signed) Mahonian polynomials over $\mathcal{D}^{\Delta}_n$, which  will be  accomplished in the following section.

\section{The (signed) Mahonian polynomials over $\mathcal{D}^{\Delta}_n$}

 This section is dedicated to deriving the formula for the (signed) Mahonian polynomials over $\mathcal{D}^{\Delta}_n$. The proofs of Theorem \ref{ThDn+} and Theorem \ref{ThDn-} are subsequently established by applying Theorems \ref{dBn}, \ref{ThBn}, and \ref{forMahddelta} into Lemma \ref{relBD-thm}.

\begin{thm} \label{forMahddelta} Let $\varepsilon=\pm1$ and for  $n\geq 2$,
\begin{align}\label{forMahddelta-eq}
   d^{\Delta, \varepsilon}_{n} (q)&= \sum_{\sigma\in \mathcal{D}^{\Delta}_{n}}\varepsilon ^ { \ell _ { B } ( \sigma ) } q^{\operatorname { fmaj }(\sigma)}  \\[5pt]
  &=[2]_q[4]_q\cdots[2n-2]_q\sum^{n-2}_{k=0}\frac{(-1)^kq^{k^2+k+1}[n-k-1]_q}{[2]_q[4]_q\cdots[2k]_q}\varepsilon^{n-k}   \left(\frac{1-q}{1-\varepsilon q}\right)^{\lceil \frac{n-k}{2} \rceil} \nonumber
\end{align}
with the convention that $ d^{\Delta, \varepsilon}_1(q)=0$.
\end{thm}
We will follow the Wachs' approach to provide a proof of  Theorem \ref{forMahddelta}. To begin, we examine the set
\[\Delta^{<} _ { n }=\{\gamma=\gamma_1\cdots \gamma_n \in B _ { n }, 0<\gamma_n<n\}.\]
We next establish the formula for the (signed) Mahonian polynomial over the set $\Delta^{<} _ { n }$.
\begin{prop}\label{D-gen}
Let $\varepsilon=\pm1$ and for $n\geq 2$,
\begin{align}\label{deltapo}
\Delta^{<,\varepsilon} _ { n }(q)=\sum _ {\gamma \in  \Delta^{<} _ { n }} \varepsilon^{\ell_B(\gamma)} q ^ { \operatorname {fmaj } ( \gamma ) } &=\varepsilon^n  q [2]_{\varepsilon q}[4]_q\cdots[2n-2]_{\varepsilon^{n-1} q}[n-1]_{\varepsilon^{n} q}   \nonumber \\[5pt]
&=\varepsilon^n  q \left(\frac{1-q}{1-\varepsilon q}\right)^{\lceil \frac{n}{2} \rceil}[2]_{q}[4]_q\cdots[2n-2]_{ q}[n-1]_{ q}.
\end{align}
with the convention that $\Delta^{<,\varepsilon} _ { 1 }(q)=0$.
\end{prop}
\proof By definition, we have
\begin{align}\label{deltar}
    \sum _ { \gamma \in \Delta^{<} _ { n } } \varepsilon^{\ell_B(\gamma)} q ^ { \operatorname { fmaj } ( \gamma ) } &=\sum _ { \gamma \in \Delta _ { n } } \varepsilon^{\ell_B(\gamma)} q ^ { \operatorname { fmaj} ( \gamma ) }-\sum _ { \gamma \in \Delta _ { n }  \atop \gamma_n=n} \varepsilon^{\ell_B(\gamma)} q ^ { \operatorname { fmaj } ( \gamma ) }.
\end{align}
Combining \eqref{delta} and \eqref{deltas}, we have
\begin{equation}\label{deltasc}
\sum _ { \gamma \in \Delta _ { n } } \varepsilon^{\ell_B(\gamma)} q ^ { \operatorname { fmaj } ( \gamma ) }=[2]_{\varepsilon q}[4]_q\cdots[2n-2]_{\varepsilon^{n-1}q}[n]_{\varepsilon^nq}.
\end{equation}
It is easy to see that
\[\sum _ { \gamma \in \Delta _ { n }  \atop \gamma_n=n} \varepsilon^{\ell_B(\gamma)} q ^ { \operatorname { fmaj } ( \gamma ) }=\sum _ { \gamma' \in B _ { n-1 }} \varepsilon^{\ell_B(\gamma')} q ^ { \operatorname { fmaj } ( \gamma' ) }.\]
From Theorem \ref{BM-1} and Theorem \ref{BMS-1}, we derive that
\begin{align}\label{deltasplu}
\sum _ { \gamma \in \Delta _ { n }  \atop \gamma_n=n} \varepsilon^{\ell_B(\gamma)} q ^ { \operatorname { fmaj } ( \gamma ) }=  [2]_{\varepsilon q}[4]_q\cdots[2n-2]_{\varepsilon^{n-1} q}.
\end{align}
Applying \eqref{deltasc} and \eqref{deltasplu} to \eqref{deltar}, we arrive at  \eqref{deltapo}. This completes the proof. \qed

We proceed to demonstrate the following proposition.
\begin{prop} \label{insertlemma-d}
   Let  $\varepsilon=\pm1$, $1\leq k\leq n$ and $\sigma\in D^{\Delta}_k$, we have
    \begin{equation} \label{insertlemma-d-eq}
    \sum_{dp(\pi)=\sigma  \atop \pi \in \Delta^{<} _ { n }} \varepsilon^{\ell_B(\pi)}q^{\operatorname { fmaj } ( \gamma )}=\varepsilon^{\ell_B(\sigma)}q^{2(n-k)+{\operatorname { fmaj } ( \sigma )}} {n-1\brack k-1}_{q^2}  .
\end{equation}
\end{prop}
The proof of Proposition \ref{insertlemma-d} is analogous to Proposition \ref{DBn} given by Chow \cite{Chow-2006}, which generalized Wachs' proof \cite{Wachs-1989}.  Let us review
the combinatorial settings of Chow.

 For $\sigma = \sigma_1  \ldots \sigma_n \in B_n$, we say that a letter $\sigma_i$ of $\sigma$ is an excedant (resp. subcedant) of $\sigma$ if $\sigma_i > i$ (resp. $\sigma_i < i$). Let $s(\sigma)$ and $e(\sigma)$ denote the numbers of subcedants and excedants of $\sigma$ respectively. It is clear that excedants of $\sigma$ are necessarily positive. 
 
 We now fix $n$ and let $k \leq n$.   {\it The map $\overline{\psi}_{n,k}$} is defined as follows: Let $\sigma \in B_k$,   $\overline{\psi}_{n,k}(\sigma)=\tilde{\sigma}$ is obtained from $\sigma$ by replacing its $i$th smallest (in absolute value) subcedant $\sigma_j$ by $(\text{sgn}(\sigma_j))i$, $i = 1, 2, \ldots, s(\sigma)$, its $i$th smallest fixed point by $s(\sigma) + i$, $i = 1, 2, \ldots, k - s(\sigma) - e(\sigma)$, and its $i$th largest excedant by $n - i + 1$, $i = 1, 2, \ldots, e(\sigma)$.

  For example, let $\sigma = 1\, 6\, \bar{3}\, 5\, 8\, 2\, 7\, \bar{4} \in B_8$ and $n=9$. We see that $s(\sigma)=3$, $e(\sigma)=3$, and 
 \[\overline{\psi}_{9,8}({\sigma})=4\, 8\, \bar{2}\, 7\, 9\, 1\, 5\, \bar{3}.\]

Chow \cite{Chow-2006} demonstrated that the map  $\overline{\psi}_{n,k}$ is the bijection between the following two sets:
\begin{itemize}
\item  For $\sigma \in {\mathcal{D}}^B_k$, let $\mathfrak{S}_n(\sigma)$ denote the set of signed permutations $\pi$ in $B_n$ such that $dp(\pi)=\sigma$.
For example,
\begin{align*}
    \mathfrak{S}_4(\overline{2}\,\overline{1})=\{1\,2\,\overline{4}\,\overline{3}, 1\,\overline{4}\,3\,\overline{2}, 1\,\overline{3}\,\overline{2}\,4, \overline{4}\,2\,3\,\overline{1}, \overline{3}\,2\,\overline{1}\,4, \overline{2}\,\overline{1}\,3\,4 \}.
\end{align*}
\item For two disjoint signed permutations $\sigma=\sigma_1\cdots\sigma_{m}$ and $\pi=\pi_1\cdots \pi_n$, we say that $\alpha$ is a shuffle of $\sigma$ and $\pi$ if both $\sigma$ and $\pi$  are
subsequences of $\alpha$. The set of shuffles of  $\sigma$ and $\pi$  is denoted $\mathfrak{S}(\sigma, \pi)$. For example,
\begin{align*}
    \mathfrak{S}(\overline{4}\,\overline{1},2\,3)=\{2\,3\,\overline{4}\,\overline{1}, 2\,\overline{4}\,3\,\overline{1}, 2\,\overline{4}\,\overline{1}\,3, \overline{4}\,2\,3\,\overline{1}, \overline{4}\,2\,\overline{1}\,3, \overline{4}\,\overline{1}\,2\,3\}.
\end{align*}
\end{itemize}

\noindent{{\bf Remark.} We say that two signed permutations $\sigma=\sigma_1\cdots\sigma_{m}$ and $\pi=\pi_1\cdots \pi_n$ are disjoint if the permutations $|\sigma_1|\cdots|\sigma_{m}|$ and $|\pi_1|\cdots |\pi_n|$ are disjoint.}

\begin{lem}[Chow]\label{Chowshuffle} For $0\leq k\leq n$, let  $\sigma \in {\mathcal{D}}^B_k$, $\tilde{\sigma}=\overline{\psi}_{n,k}(\sigma)$, and $\gamma=s(\sigma)+1,s(\sigma)+2,\ldots,n-e(\sigma)$. Then the map $\overline{\psi}_{n,n}$ is a bijection between $\mathfrak{S}_n(\sigma)$ and $\mathfrak{S}(\tilde{\sigma},\gamma)$. Moreover, for $\pi \in \mathfrak{S}_n(\sigma)$, we have $\overline{\psi}_{n,n}(\pi) \in \mathfrak{S}(\tilde{\sigma},\gamma)$ such that
 ${\rm fmaj}(\pi)={\rm fmaj}(\overline{\psi}_{n,n}(\pi)).$
\end{lem}
Note that  a generalized version of Lemma \ref{Chowshuffle} for the wreath product $\mathcal{C}_r \wr\mathfrak{S}_n$ can be found in \cite{Assaf-2010}.

Subsequently, Chow \cite{Chow-2006} applied the following formula due to Garsia and Gessel  \cite{Garsia-Gessel-1979} in the context of  the hyperochaherdral group $B_n$. Specifically, for two disjoint permutations $\sigma=\sigma_1\cdots\sigma_{m}$ and $\pi=\pi_1\cdots \pi_n\in $, Garsia and Gessel \cite{Garsia-Gessel-1979} showed that
\begin{equation}\label{GG-F}
    \sum_{\alpha \in \mathfrak{S}(\sigma,\pi)}q^{{\rm maj}(\alpha)}=
   q^{{\rm maj}(\sigma)+{\rm maj}(\pi)}{n+m \brack n}_q.
\end{equation}
Note that Garsia and Gessel's formula \eqref{GG-F} can also be derived from Stanley's shuffle theorem \cite{Stanley-1972} by employing $q$-analogue of the Chu-Vandermonde summation (see \cite[Eq.3.3.10]{Andrews-1976}),
\begin{equation}\label{qchuvan}
\sum_{k=0}^h {n \brack k}_q{m \brack h-k}_q q^{(n-k)(h-k)}={m+n \brack h}_q.
\end{equation}
For more information, we refer to \cite{Ji-Zhang-2024}.

 Chow \cite{Chow-2006} extended the formula  \eqref{GG-F} into the following form:
\begin{lem}[Garsia-Gessel-Chow]\label{Chowshufflebb}For two disjoint signed permutations $\sigma=\sigma_1\cdots\sigma_{m}$ and $\pi=\pi_1\cdots \pi_n\in $, we have
\begin{equation}
    \sum_{\alpha \in \mathfrak{S}(\sigma,\pi)}q^{{\rm fmaj}(\alpha)}=
   q^{{\rm fmaj}(\sigma)+{\rm fmaj}(\pi)}{n+m \brack n}_{q^2}.
\end{equation}
\end{lem}
By combining Lemma \ref{Chowshuffle} and Lemma \ref{Chowshufflebb}, Chow \cite{Chow-2006} successfully concluded the proof of    Proposition \ref{DBn}.

In order to justify Proposition \ref{insertlemma-d}, we first show that the map  $\overline{\psi}_{n,k}$ can be refined into the following two subsets.
\begin{itemize}
\item For $\sigma \in {\mathcal{D}}^\Delta_k$, let ${\Delta}^{<} _ { n }(\sigma)$ denote the set of signed permutations $\pi=\pi_1\cdots \pi_n$ in $B_n$ such that $dp(\pi)=\sigma$ and $0<\pi_n<n$.
For example,
\begin{align*}
   \Delta^{<} _ { 5 }(2\,\overline{3}\,1)=\{1\,2\,4\,\overline{5}\,3,1\,4\,3\,\overline{5}\,2,1\,3\,\overline{5}\,4\,2,4\,2\,3\,\overline{5}\,1,3\,2\,\overline{5}\,4\,1,2\,\overline{5}\,3\,4\,1\}.
\end{align*}

\item For two disjoint signed permutations $\sigma=\sigma_1\cdots\sigma_{m}$ and $\pi=\pi_1\cdots \pi_n$, let
$\mathfrak{S}^{sb}(\sigma,\pi)$ denote the set of shuffles
$\alpha=\alpha_1\cdots\alpha_{n+m}$ of $\sigma$ and
$\pi$ such that $\alpha_{n+m}=\min \{\pi_n,\sigma_m\}$. As an illustration, we see that
\begin{align*}
    \mathfrak{S}^{sb}(4\,\overline{5}\,1,2\,3)=\{2\,3\,4\,\overline{5}\,1,2\,4\,3\,\overline{5}\,1,2\,4\,\overline{5}\,3\,1,4\,2\,3\,\overline{5}\,1,4\,2\,\overline{5}\,3\,1,4\,\overline{5}\,2\,3\,1\}.
\end{align*}
\end{itemize}

\begin{lem}\label{Chowshuffleref}
 For $1\leq k\leq n$, let  $\sigma \in {\mathcal{D}}^\Delta_k$, $\tilde{\sigma}=\overline{\psi}_{n,k}(\sigma)$, and $\gamma=s(\sigma)+1,s(\sigma)+2,\ldots,n-e(\sigma)$. Then the map $\overline{\psi}_{n,n}$ is a bijection between $\Delta^{<}_n(\sigma)$ and $\mathfrak{S}^{sb}(\tilde{\sigma},\gamma)$. Moreover, for $\pi \in \Delta^{<}_n(\sigma)$, we have $\overline{\psi}_{n,n}(\pi) \in \mathfrak{S}^{sb}(\tilde{\sigma},\gamma)$ such that ${\rm fmaj}(\pi)={\rm fmaj}(\overline{\psi}_{n,n}(\pi)).$
 \end{lem}

\pf
Let $\sigma=\sigma_1\cdots \sigma_k$ be a signed derangement in $\mathcal{D}^{\Delta}_k$. By definition, we see that $0<\sigma_k<k$ is a subcedant. Let $\tilde{\sigma}=\overline{\psi}_{n,k}(\sigma)$, from the construction of $\overline{\psi}_{n,k}$, we derive that  $\tilde{\sigma}=\tilde{\sigma}_1\cdots \tilde{\sigma}_k$ is a signed permutation of $\{1, 2, \ldots, s(\sigma)\} \cup \{n - e(\sigma) + 1, n - e(\sigma) + 2, \ldots, n\}$ such that $\tilde{\sigma}_k\le s(\sigma)$. It follows that $0<\tilde{\sigma}_k<\gamma_{n-k}.$

Let $\pi=\pi_1\cdots \pi_n \in \Delta^{<}_n(\sigma)$. By definition, we see that $\pi_n<n$. Let $\tilde{\pi}=\overline{\psi}_{n,n}(\pi)$.  From Lemma \ref{Chowshuffle}, we   know  that $\tilde{\pi}=\tilde{\pi}_1\cdots \tilde{\pi}_n \in \mathfrak{S}(\tilde{\sigma},\gamma)$ and ${\mathrm{fmaj}}(\pi)={\mathrm{fmaj}}(\tilde{\pi})$. Since $\pi_n<n$, and from the construction of $\overline{\psi}_{n,n}$, we see that $\tilde{\pi}_n\le s(\pi)=s(\sigma)$. Hence we derive that $\tilde{\pi}_n=\tilde{\sigma}_k$, and so $\tilde{\pi} \in \mathfrak{S}^{sb}(\tilde{\sigma},\gamma)$.

Conversely, let $\tilde{\pi}=\tilde{\pi}_1\cdots \tilde{\pi}_n  \in \mathfrak{S}^{sb}(\tilde{\sigma},\gamma)$. In this case, we see that $\tilde{\pi}_n=\tilde{\sigma}_k$. Let $\pi=\overline{\psi}^{-1}_{n,n}(\tilde{\pi})$. Using Lemma \ref{Chowshuffle}, we derive that $\pi \in B_n$ and $dp(\pi)=\sigma$.  Since $\tilde{\pi}_n=\tilde{\sigma}_k>0$, by  the construction of $\overline{\psi}_{n,n}$, we derive that $0<\pi_n<n$. Consequently, $\pi\in \Delta^{<}_n(\sigma)$.  This completes the proof. \qed

Next, we establish the following refinement of   Garsia and Gessel's formula in the hyperochaherdral group $B_n$, which can be derived from a refinement of Stanley's shuffle theorem established in \cite{Ji-Zhang-2024}.
\begin{lem}\label{gg-r-thm}  For two disjoint signed permutations $\sigma=\sigma_1\cdots\sigma_{m}$ and $\pi=\pi_1\cdots \pi_n$ with   $\sigma_m\prec \pi_n$ in the order given by \eqref{neword}, we have
\begin{equation}\label{gg-r-e}
    \sum_{\alpha \in \mathfrak{S}^{sb}(\sigma,\pi)}q^{{\rm fmaj}(\alpha)}=
  q^{{\rm fmaj}(\sigma)+{\rm fmaj}(\pi)+2n} {n+m-1 \brack n}_{q^2}.
\end{equation}
\end{lem}
\pf Let us first recall the refinement of Stanley's shuffle theorem established in \cite{Ji-Zhang-2024}. Assume that $\sigma=\sigma_1\cdots\sigma_{m}\in  \mathfrak{S}_{A}$ and $\pi=\pi_1\cdots \pi_n\in  \mathfrak{S}_{B}$ be two disjoint permutations with $\pi_n>\sigma_m$. We have
 \begin{align}\label{stan-r}
    \sum_{\alpha\in \mathfrak{S}^{sb}(\sigma,\pi) \atop {\rm des}(\alpha)=k}q^{{\rm maj}(\alpha)}&=
    {m-{\rm des}(\sigma)+{\rm des}(\pi) \brack k-{\rm des}(\sigma)}_q {n-{\rm des}(\pi)+{\rm des}(\sigma)-1 \brack n-k+{\rm des}(\sigma)}_q \nonumber\\[5pt]
    &\quad \quad \quad \times q^{{\rm maj}(\sigma)+{\rm maj}(\pi)+n+(k-{\rm des}(\pi)-1)(k-{\rm des}(\sigma))}.
\end{align}
By employing $q$-analogue of the Chu-Vandermonde summation \eqref{qchuvan}, we derive from \eqref{stan-r} that
 \begin{align}\label{GG-r}
    \sum_{\alpha\in \mathfrak{S}^{sb}(\sigma,\pi)}q^{{\rm maj}(\alpha)}=
    q^{{\rm maj}(\sigma)+{\rm maj}(\pi)+n}{m+n-1 \brack n}_q.
\end{align}
We proceed to generalize this result to the shuffle of two signed permutations. To distinguish, here we let $\overline{\sigma}=\overline{\sigma}_1\cdots \overline{\sigma}_m \in \mathfrak{S}_{A}$ and   $\overline{\pi}=\overline{\pi}_1\cdots\overline{\pi}_n \in \mathfrak{S}_{B}$ be  two disjoint signed permutations. Assume that $\overline{\alpha}=\overline{\alpha}_1\cdots \overline{\alpha}_{n+m}$ is the shuffle of $\overline{\sigma}$ and $\overline{\pi}$. From the definition of the shuffle, it is easy to see that
\begin{equation}\label{rel-neg}
{\rm neg}(\overline{\alpha})={\rm neg}(\overline{\sigma})+{\rm neg}(\overline{\pi}).
\end{equation}
Assume that $\overline{\sigma}_m\prec \overline{\pi}_n$ in the order given by \eqref{neword}, it is not difficult to derive from  \eqref{GG-r} and \eqref{rel-neg} that
\begin{align}\label{pf-prop}
    \sum_{\overline{\alpha}\in \mathfrak{S}^{sb}(\overline{\sigma},\overline{\pi})}z^{{\rm neg}(\overline{\alpha})}q^{{\rm maj}(\overline{\alpha})}=
    z^{{\rm neg}(\overline{\sigma})+{\rm neg}(\overline{\pi})}q^{{\rm maj}(\overline{\sigma})+{\rm maj}(\overline{\pi})+n}{m+n-1 \brack n}_q,
\end{align}
where   ${\rm maj}(\overline{\alpha})$ is the major index of $\overline{\alpha}$ in the definition \eqref{defi-fmaj} of ${\rm fmaj}$, computed by using the order given by  \eqref{neword}.
 By substituting $q$ with $q^2$ in \eqref{pf-prop} and setting $z=q$,  according to the definition \eqref{defi-fmaj} of ${\rm fmaj}$, we derive \eqref{gg-r-e}. This completes the proof. \qed

With Lemma \ref{Chowshuffleref} and Lemma \ref{gg-r-thm} at our disposal, we are prepared to establish Proposition \ref{insertlemma-d}. As a result, we can prove Theorem \ref{forMahddelta}.

 \noindent{\it Proof of Proposition \ref{insertlemma-d}.} Utilizing Lemma \ref{Chowshuffleref} and applying \eqref{bnmmtt}, we derive that for $\sigma \in \mathcal{D}^{\Delta}_k$, 
 \begin{align}\label{insertlemma-d-eq-pf}
    \sum_{dp(\pi)=\sigma  \atop \pi \in \Delta^{<} _ { n }} \varepsilon^{\ell_B(\pi)}q^{\operatorname { fmaj } ( \gamma )}&=\varepsilon^{\ell_B(\sigma)} \sum_{\alpha \in \mathfrak{S}^{sb}(\tilde{\sigma},\gamma)}q^{{\rm fmaj}(\alpha)},
\end{align}
 It is noteworthy that $\tilde{\sigma}$ has $k$ elements, $\gamma$ has $n-k$ elements, and ${\rm fmaj}(\gamma)=0$. By incorporating Lemma \ref{gg-r-thm} into \eqref{insertlemma-d-eq-pf}, we consequently establish \eqref{insertlemma-d-eq}. This completes the proof. \qed

\noindent{\it Proof of Theorem \ref{forMahddelta}.} By summing over all signed derangements $\sigma \in \mathcal{D}^{\Delta}_k$ for $1 \leq k \leq n$, we infer from Proposition \ref{insertlemma-d} that
\begin{align*}
&\sum_{k=1}^n{n-1\brack k-1}_{q^2} q^{2(n-k)} \sum_{\sigma \in \mathcal{D}^{\Delta}_k} \varepsilon^{\ell_B(\sigma)}q^{\operatorname { fmaj } ( \sigma )}\\[5pt]
&=\sum_{k=1}^n\sum_{\sigma \in \mathcal{D}^{\Delta}_k}\sum_{ \gamma \in \Delta^{<} _ { n }(\sigma) } \varepsilon^{\ell_B(\gamma)}q^{\operatorname { fmaj } ( \gamma )}\\[5pt]
&= \sum_{  \gamma \in \Delta^{<} _ { n } } \varepsilon^{\ell_B(\gamma)}q^{\operatorname { fmaj } ( \gamma )},
\end{align*}
which can be further simplified as
\begin{align*}
\sum_{k=0}^{n}{n\brack k}_{q^2} q^{-2k} \sum_{\sigma \in \mathcal{D}^{\Delta}_{k+1}} \varepsilon^{\ell_B(\sigma)}q^{\operatorname { fmaj } ( \sigma )}&= q^{-2n}\sum_{  \gamma \in \Delta^{<} _ { n+1 } } \varepsilon^{\ell_B(\gamma)}q^{\operatorname { fmaj } ( \gamma )}.
\end{align*}
By employing $q$-binomial inversion \cite[Corollary 3.38]{Aigner-1979}, we deduce  that for $n\geq 0$,
\begin{align} \label{d-inv-pfa}
q^{-2n} \sum_{\sigma \in \mathcal{D}^{\Delta}_{n+1}} \varepsilon^{\ell_B(\sigma)}q^{\operatorname { fmaj } ( \sigma )}&= \sum_{k=0}^{n}(-1)^k{n\brack k}_{q^2} q^{2{k\choose 2}} q^{-2(n-k)}\sum_{  \gamma \in \Delta^{<} _ { n-k+1 } } \varepsilon^{\ell_B(\gamma)}q^{\operatorname { fmaj } ( \gamma )}.
\end{align}
Applying \eqref{deltapo} to \eqref{d-inv-pfa}, we obtain \eqref{forMahddelta-eq}.
This completes the proof. \qed

We conclude this paper with the proofs of Theorem \ref{ThDn+} and Theorem \ref{ThDn-} achieved through the utilization of Theorems \ref{dBn}, \ref{ThBn}, and \ref{forMahddelta} in conjunction with Lemma \ref{relBD-thm}.

\noindent{\it Proofs of Theorem \ref{ThDn+} and Theorem \ref{ThDn-}:}  From Theorem \ref{forMahddelta}, we have
\begin{align}\label{neq-n}
    &\sum_{\sigma\in \mathcal{D}^{\Delta}_n}\varepsilon^ { \ell _ { B } ( \sigma ) } (\varepsilon q)^{\operatorname { fmaj }(\sigma)}\notag\\[5pt]
    &=[2]_{\varepsilon q}[4]_{\varepsilon q}\cdots[2n-2]_{\varepsilon q}\sum^{n-2}_{k=0}\frac{(-1)^k(\varepsilon q)^{k^2+k+1}[n-k-1]_{\varepsilon q}}{[2]_{\varepsilon q}[4]_{\varepsilon q}\cdots[2k]_{\varepsilon q}}\varepsilon^{n-k}   \left(\frac{1-\varepsilon q}{1-q}\right)^{\lceil \frac{n-k}{2} \rceil}\notag\\[5pt]
    &=[2]_q[4]_q\cdots[2n-2]_q\sum^{n-2}_{k=0}\frac{(-1)^kq^{k^2+k+1}(\varepsilon^{n-k-1}-q^{n-k-1})}{[2]_q[4]_q\cdots[2k]_q(1-q)} \left(\frac{1-q}{1-\varepsilon q}\right)^{\lfloor \frac{n-k}{2}\rfloor}.
\end{align}
Combining  Theorem \ref{dBn} and Theorem \ref{ThBn} leads to
\begin{align}\label{=-n}
&\frac{\varepsilon}{2}\left(\sum_{\pi\in \mathcal{D}^B_{n-1} }\varepsilon^{\ell_B(\pi)}q^{\operatorname { fmaj }(\pi)}-\sum_{\pi\in \mathcal{D}^B_{n-1} }\varepsilon^{\ell_B(\pi)}(-q)^{\operatorname { fmaj }(\pi)}\right)\nonumber \\[5pt]
    &=[2]_q[4]_q\cdots[2n-2]_q\sum^{n-1}_{k=0}\frac{(-1)^kq^{2\binom{k}{2}}}{[2]_q[4]_q\cdots[2k]_q}\notag\\[5pt]
    &\times \frac{\varepsilon}{2}\left(\left(\frac{1-q} {1-\varepsilon q}\right)^{\lfloor \frac{n-k}{2}\rfloor}-\left(\frac{1-q} {1+ q}\right)^{ n-1-k}\left(\frac{1+q} {1+\varepsilon q}\right)^{\lfloor \frac{n-k}{2}\rfloor}\right).
\end{align}
Upon substituting \eqref{neq-n} and \eqref{=-n} with Lemma \ref{relBD-thm} for $\varepsilon=1$,
 we arrive at
\begin{align*}
    \sum_{\pi\in \mathcal{D}^{D}_n}q^{\operatorname { Dmaj }(\pi)}&=[2]_q[4]_q\cdots[2n-2]_q\sum^{n-2}_{k=0}\frac{(-1)^kq^{k^2+k+1}[n-1-k]_q}{[2]_q[4]_q\cdots[2k]_q}\\
    &+[2]_q[4]_q\cdots[2n-2]_q\sum^{n-1}_{k=0}\frac{(-1)^kq^{2\binom{k}{2}}}{[2]_q[4]_q\cdots[2k]_q}\frac{1}{2}\left(1-\left(\frac{1-q} {1+ q}\right)^{ n-k-1}\right)\\
    &=[2]_q[4]_q\cdots[2n-2]_q\sum^{n-2}_{k=0}\frac{(-1)^kq^{k^2-k}}{[2]_q[4]_q\cdots[2k]_q} \\
    &\times \left(q^{2k+1}[n-1-k]_q+ \frac{1}{2}\left(1-\left(\frac{1-q} {1+ q}\right)^{ n-k-1}\right) \right),
\end{align*}
which is in accordance with \eqref{ThDn+-eq}.

On the other hand,  plugging  \eqref{neq-n} and \eqref{=-n} into Lemma \ref{relBD-thm} for $\varepsilon=-1$, we derive that
\begin{align*}
    &\sum_{\pi\in \mathcal{D}^{D}_n}(-1)^{\ell_D(\pi)}q^{\operatorname { Dmaj }(\pi)}\\[5pt]
        &=[2]_q[4]_q\cdots[2n-2]_q\sum^{n-2}_{k=0}\frac{(-1)^{n-1}q^{2\binom{k}{2}}}{[2]_q[4]_q\cdots[2k]_q} \left(\frac{1-q}{1+ q}\right)^{\lfloor \frac{n-k}{2}\rfloor}\\[5pt]
     &\quad \times\left( q^{2k+1}\frac{(1+(-1)^{n-k}q^{n-k-1})}{1-q}+\frac{(-1)^{n-k}}{2}\left(1-\left(\frac{1-q} {1+ q}\right)^{ n-1-k}\left(\frac{1+q} {1- q}\right)^{2\lfloor \frac{n-k}{2}\rfloor}\right)\right)\\[5pt]
     &=[2]_q[4]_q\cdots[2n-2]_q\sum^{n-2}_{k=0}\frac{(-1)^{n-1}q^{2\binom{k}{2}}}{[2]_q[4]_q\cdots[2k]_q} \left(\frac{1-q}{1+ q}\right)^{\lfloor \frac{n-k}{2}\rfloor}\\[5pt]
     &\quad\times\frac{2q^{2k+1}(1+(-1)^{n-k}q^{n-k-1})-(1+(-1)^{n-k})q}{2(1-q)},
\end{align*}
 aligning with \eqref{ThDn--eq}. This completes the proof. \qed

 \vskip 0.2cm
\noindent{\bf Acknowledgment.} We are grateful to the anonymous referees for their insightful comments and suggestions. This work
was supported by   the National   Science Foundation of China.

\end{document}